\documentclass[12pt]{amsart}
\usepackage{a4wide,amssymb}

\newtheorem{thm}{Theorem}
\newtheorem{cor}[thm]{Corollary}

\newtheorem{prop}[thm]{Proposition}
\newtheorem{defn}[thm]{Definition}

\newcommand\Hom{\operatorname{Hom}}
\newcommand\Der{\operatorname{Der}}
\newcommand\Ker{\operatorname{Ker}}
\newcommand\Coker{\operatorname{Coker}}
\newcommand\Tor{\operatorname{Tor}}
\newcommand\Ext{\operatorname{Ext}}

\def\two{\begin{matrix}\to\\[-3mm]\to\end{matrix}}

\def\three{\begin{matrix}\to\\[-2mm]\to\\[-2mm]\to\end{matrix}}

\def\dtwo{\begin{matrix}\to\\[-3.5mm]\vdots\\[-3mm]\to\end{matrix}}

\title{Equivariant homology and cohomology of groups }%
\author{Hvedri Inassaridze}%
\address{A.Razmadze Mathematical Institute, Georgian Academy of Sciences,
 M.Alexidze Str.1, Tbilisi 0193, Georgia}%
\email{hvedri@rmi.acnet.ge}%


\advance\headheight1.15pt

\begin{document}
\maketitle

\

\begin{abstract} We provide and study an equivariant theory of group (co)homology
of a group $G$ with coefficients in a $\Gamma$-equivariant
$G$-module $A$, when a separate group $\Gamma$ acts on $G$ and
$A$, generalizing the classical Eilenberg-MacLane (co)homology
theory of groups. Relationship with equivariant cohomology of
topological spaces is established and application to algebraic
K-theory is given.
\end{abstract}

\

\

\section*{Introduction}
It is well known that the study of groups with operators has many
important applications in algebra and topology. The category of
groups enriched with an action by automorphisms of a given group
provides a suitable setting for the investigation of an extensive
list of subjects with recognized mathematical interest. See for
instance recent results in equivariant stable homotopy theory
\cite{Ca} and articles devoted to equivariant algebraic K-theory
\cite{FHM, Ph, Ku}. The origin of the equivariant investigation in
homological algebra, particularly in extension theory of groups,
goes back to the article of J.H.C.Whitehead \cite{Wh}. It should
be noted that recently a theory of cohomology of groups with
operators was developed \cite{CGO2}, motivated by the graded
categorical groups classification problem which was suggested by
Frohlich and Wall \cite{FrhWa}. This problem was solved
\cite{CGO1} by using the third cohomology of groups with operators
introduced in \cite{CGO2}. An equivariant version of the classical
Brauer-Hasse-Noether result was proved \cite{CeGa1} showing that
for any Galois finite field extension $F/K$ on which a separate
group of operators $\Gamma$ is acting, there is an isomorphism of
equivariant isomorphism classes of finite dimensional central
simple $K$-algebras endowed with a $\Gamma$-action and containing
$F$ as an equivariant strictly maximal subfield and the second
cohomology of groups with operators defined in \cite{CGO2} of the
Galois group of the extension. A homology theory of groups with
operators corresponding to the cohomology theory of groups with
operators \cite{CGO2} has been treated in \cite{CeIn}.

In \cite{CGO2} it was stated that the second cohomology group of a
$\Gamma$-group $G$ with coefficients in a $\Gamma$-equivariant
$G$-module $A$ classifies the $\Gamma$-equivariant extensions of
$G$ by $A$. From this result arises the natural problem about the
cohomological characterization of those $\Gamma$-equivariant
extensions of $G$ by $A$ which are $\Gamma$-splitting. The
solution of this problem (see Theorem \ref{t-20}) has motivated an
attempt to develop a different equivariant (co)homology theory of
groups, which is presented in this paper.

By definition for a $\Gamma$-group $G$ its equivariant homology
and cohomology groups, $H_{n}^{\Gamma}(G,-)$ and
$H^{n}_{\Gamma}(G,-)$, are defined as relative
$\Tor^{\mathcal{F}}_{n}$ and $\Ext^{n}_{\mathcal{F}}$, $n\geq 0$,
functors respectively in the category of $\Gamma$-equivariant
$G$-modules (Definition \ref{d-1}). Therefore this (co)homology
theory of groups can be considered as a part of the relative
homological algebra \cite{El}. We provide equivariant versions of
classical homological theorems: (co)chain and cotriple
presentations of the homology and cohomology of groups, Hopf
formula for the second integral homology, universal coefficient
formulas, universal central extensions, cohomological
classification of extensions of groups, exact (co)homology
sequences, Tate (co)homology of groups and the cup product.
Applications in algebraic K-theory (Corollary 24) and the
relationship with equivariant cohomology of topological spaces
(Theorem 22) are established.

Corollary 24 motivates the following

\medskip

\noindent{\bf Conjecture.} {\em There is an isomorphism
$K_{3}(A)\cong H^{St(\mathbb{Z})}_{3}(St(A))$ for any ring $A$.}

\medskip

For its proof equivariant versions of relevant classical homotopy
theorems will be probably needed. This is well known when $A$ is a
unital ring \cite{Ge} and in this case $St(\mathbb{Z})$ acts
trivially on $St(A)$.

Another application will be the construction of an alternative
equivariant algebraic K-theory $K^{\Gamma}_{*}$ by using
$\Gamma$-equivariant commutators ( Section 6). Moreover in the
near future it is intended to investigated for any ring (not
necessarily with unit) the relationship of higher Quillen's
algebraic K-groups with the equivariant integral homology of the
general linear group under the action of the Steinberg group of
the ring of integers. It is also intended to establish the
relationship of this alternative equivariant algebraic K-theory
with existing equivariant algebraic K-theory and equivariant
homotopy theory and to provide higher Hopf formulae for
equivariant integral homology of groups.

\

\

\section{Definition and (co)cycle description of the equivariant
(co)homology of groups}

\

Before defining the equivariant (co)homology of groups we briefly
recall the definition of the cohomology of groups with operators
introduced in \cite{CGO2}.

The $n$-cochains of a $\Gamma$-group $G$ with coefficients in a
$\Gamma$-equivariant $G$-module $A$ (see definition below) are the
maps
$$
 f: \bigcup G^{p+1}\times \Gamma^{q}\rightarrow A, \;\; p+q=
n-1,
$$
normalized in the sense that $f(x_{1},\cdots,
x_{p+1},\sigma_{1},\cdots,\sigma_{q})=0$ whenever $x_{i}=1$ or
$\sigma_{j}=1$ for some  $i=1,\ldots,p+1$ or $j=1,\ldots,q$, and
the coboundary operator is introduced in a natural way. One gets a
cochain complex whose homology groups are the cohomology groups of
the $\Gamma$-group $G$. Note that in this theory the zero
cohomology group is trivial and as mentioned above the second
cohomology group describes the $\Gamma$- equivariant extensions of
the $\Gamma$-group $G$ by the $G\rtimes\Gamma$-module $A$.

Now we will give our definition of the equivariant (co)homology of
groups. Let $G$ be a $\Gamma$-group. A $\Gamma$-equivariant
$G$-module $A$ is a $G$-module equipped with a $\Gamma$-module
structure and both actions of $G$ and $\Gamma$ on $A$ satisfy the
following condition:
\begin{align}\label{1}
 ^{\sigma}(^{x}a)\; =\;  ^{^{\sigma}x}{(^{\sigma}a)},\;\;  x\in G,
\sigma\in\Gamma,\;\; a\in A .
\end{align}
The category of $\Gamma$-equivariant $G$-modules is equivalent to
the category of $G\rtimes\Gamma$-modules, where $G\rtimes \Gamma$
denotes the semidirect product of $G$ and $\Gamma$ (see
\cite{CGO2}). Let $B$ and $C$ be two $\Gamma$-equivariant
$G$-modules. Clearly a map $f:B\rightarrow C$ is a $G\rtimes
\Gamma$-homomorphism if and only if it is compatible with the
actions of $G$ and $\Gamma$. A $G\rtimes\Gamma$-module free as a
$G$-module with basis a $\Gamma$-subset will be called a
relatively free $G\rtimes\Gamma$-module. Denote by $\mathcal{P}$
the class of $G\rtimes\Gamma$-modules which are retracts of
relatively free $G\rtimes\Gamma$-modules. The elements of
$\mathcal{P}$ will be called relatively projective $G\rtimes
\Gamma$-modules. A $G\rtimes\Gamma$-homomorphism $f:B\rightarrow
C$ of $G\rtimes \Gamma$-modules is a $\mathcal{P}$-epimorphism if
it is $\Gamma$-splitting, that is there is a $\Gamma$-map
$\gamma:C\rightarrow B$ such that $f\gamma = 1_{C}$. The group
ring $Z(G)$ is a relatively free $G\rtimes\Gamma$-module in a
natural way with the action of $\Gamma$ by
$$^{\sigma}(\underset{i}{\sum} m_{i}g_{i})=\underset{i}{\sum}
m_{i}\;^{\sigma}g_{i}\;.$$

Let $A$ be a $G\rtimes\Gamma$-module. Denote by
$I_{G\rtimes\Gamma}A$ the subgroup of $A$ generated by the
elements ${}^{(g,\sigma)}a - a = {}^{g}{({}^{\sigma}a)} - a$,
$g\in G$, $\sigma\in\Gamma$, $a\in A$, and by $A_{G\rtimes\Gamma}$
the quotient group of $A$ by $I_{G\rtimes\Gamma}A$. Then it is
easily checked that one has canonical isomorphisms
$$
\mathbb{Z}(G){\otimes}_{G\rtimes \Gamma} A \cong A_{\Gamma},\quad
\mathbb{Z}\otimes_{G\rtimes \Gamma} A \cong A_{G\rtimes
\Gamma},\quad \Hom_{G\rtimes \Gamma}(\mathbb{Z}(G),A) \cong
A^{\Gamma}.
$$
Clearly if $\Gamma$ acts trivially on $A$, then
$\mathbb{Z}(G)\otimes_{G\rtimes \Gamma} A \cong \Hom_{G\rtimes
\Gamma}(\mathbb{Z}(G),A) \cong A$.

In the category of $G\rtimes \Gamma$-modules there are sufficient
relatively projective (free) $G\rtimes \Gamma$-modules. If $A$ is
a $G\rtimes\Gamma$-module, take the free $G$-module $F(A)$
generated by $A$ and define the action of $\Gamma$ on $F(A)$ by
$$
{}^{\sigma}{(g|a|)}={}^ {\sigma}{g}|{}^{\sigma}{a}|\;,\quad g\in
G\;,\; \sigma\in \Gamma\;,\; a\in A\;.
$$
Then $F(A)$ becomes a relatively free $G\rtimes\Gamma$-module with
basis $A$ being a $\Gamma$-subset of $F(A)$ and the canonical map
$F(A)\rightarrow A$ is a $\mathcal{P}$-epimorphism, since it is
$\Gamma$-splitting by the map
$$
\gamma: A\rightarrow F(A)\;,\quad \gamma(a)= |a|\;,\quad a\in A\;.
$$

It is standard to show that one has isomorphisms
$$
H_{n}(P_{*}(A)\otimes_{G\rtimes \Gamma} B) \cong
H_{n}(A\otimes_{G\rtimes \Gamma} P_{*}(B))\;,\quad n\geq 0\;,
$$
where $P_{*}(A)$ and $P_{*}(B)$ are $\mathcal{P}$-projective
$G\rtimes\Gamma$-resolutions of $A$ and $B$ respectively.

Now we are ready to define the equivariant homology
$H^{\Gamma}_{*}(G,A)$ and cohomology $H^{*}_{\Gamma}(G,A)$ of a
$\Gamma$-group $G$ with coefficients in a $\Gamma$-equivariant
$G$-module $A$.

\medskip

\begin{defn}\label{d-1}
$H^{\Gamma}_{n}(G,A) = \Tor^{\mathcal{P}}_{n}(\mathbb{Z},A)$ and
$H^{n}_{\Gamma}(G,A) = \Ext^{n}_{\mathcal{P}}(\mathbb{Z},A)$ for
$n\geq 0$, where $G$ and $\Gamma$ act trivially on $\mathbb{Z}$.
\end{defn}

\medskip

It is clear that $H^{\Gamma}_{n}(G,A) \cong H_{n}(G,A_{\Gamma})$
and $H^{n}_{\Gamma}(G,A) \cong H^{n}(G,A^{\Gamma})$ for $n\geq 0$,
if $\Gamma$ acts trivially on $G$ and therefore this case is not
interesting from the equivariant point of view.

A short exact sequence of $G\rtimes\Gamma$-modules
\begin{align}\label{3}
0\rightarrow C_{1}\rightarrow C \overset{\beta}{\rightarrow}
C_{2}\rightarrow0
\end{align}
will be called proper if $\beta$ is $\Gamma$-splitting, i.e. there
is a $\Gamma$-map $\gamma : C_{2}\rightarrow C$ such that
$\beta\gamma= 1_{C_2}$.

Let
\begin{align}\label{4}
\cdots\rightarrow B_{n}\rightarrow\cdots\rightarrow
B_{1}\rightarrow B_{0}\rightarrow \mathbb{Z}\rightarrow 0
\end{align}
be the bar resolution of $\mathbb{Z}$, where $B_{0}=
\mathbb{Z}(G)$, and $B_{n}$, $n>0$, is the free ${\mathbb
Z}(G)$-module generated by $[g_{1},g_{2},\ldots,g_{n}], g_{i}\in
G$. Define the action of the group $\Gamma$ on the bar resolution
as follows. $\Gamma$ acts trivially on $\mathbb{Z}$, the action of
$\Gamma$ on $B_{0}$ is already defined and if $n>0$ then
$^{\sigma}(g[g_{1},g_{2},\ldots,g_{n}]) = \;
^{\sigma}g[^{\sigma}g_{1},^{\sigma}g_{2},\ldots,^{\sigma}g_{n}]$
for the action of $\Gamma$ on $B_{n}$. The well known contraction
$\gamma_{-1}:{\mathbb Z}\rightarrow B_{0}$, $\gamma_{-1}(z)= z1$,
$\gamma_n : B_n\rightarrow B_{n+1}$,
$\gamma_{n}(g[g_{1},\ldots,g_{n}])= [g,g_{1},\ldots,g_{n}]$, $n\ge
0$, is clearly a $\Gamma$-map. We deduce that under this action of
$\Gamma$ the bar resolution (\ref{4}) becomes an exact sequence of
$G\rtimes\Gamma$-modules such that each $B_{n}$ is a relatively
free $G\rtimes\Gamma$-module and the sequences
$$
0\rightarrow \Ker\partial_{n}\rightarrow B_{n}\rightarrow
Im\partial_{n}\rightarrow 0\;,\quad n>0\;
$$
and
$$
0\rightarrow \Ker\epsilon\rightarrow B_{0}\rightarrow
\mathbb{Z}\rightarrow 0\;,
$$
are proper short exact sequences of
$G\rtimes\Gamma$-modules. Therefore (\ref{4}) is a relatively free
$G\rtimes\Gamma$-resolution of $\mathbb{Z}$ which will be called
the $\Gamma$-equivariant bar resolution of $\mathbb{Z}$. It
follows that
$$
H^{\Gamma}_{n}(G,A)\cong H_{n}(B_{*}\otimes_{G\rtimes\Gamma} A),
H^{n}_{\Gamma}(G,A)\cong H_{n}(\Hom_{G\rtimes\Gamma}(B_{*},A),
n\geq 0.
$$
These isomorphisms allow us to produce an alternative description
by (co)cycles of the equivariant (co)homology of groups. The case
of the equivariant homology of groups is clear. To this end
consider the abelian group $C^{n}_{\Gamma}(G,A)$ of $\Gamma$-maps
$f: G^{n}\rightarrow A$, $n>0$, which will be called the group of
$n$-th $\Gamma$-cochains. By using the classical cobord operators
$\delta^{n}: C^{n}_{\Gamma}(G,A)\rightarrow
C^{n+1}_{\Gamma}(G,A)$, $n>0$, one gets a cochain complex
$$
0\rightarrow C^{0}_{\Gamma}(G,A)\rightarrow
C^{1}_{\Gamma}(G,A)\rightarrow
C^{2}_{\Gamma}(G,A)\rightarrow\cdots\rightarrow
C^{n}_{\Gamma}(G,A)\rightarrow\cdots ,
$$
where $C^{0}_{\Gamma}= A^{\Gamma}$, $\Ker\delta^{1}=
\Der_{\Gamma}(G,A)$ is the group of $\Gamma$-derivations, and the
homology groups of this complex give the $\Gamma$-equivariant
cohomology groups of $G$ with coefficients in the $G\rtimes
\Gamma$-module $A$.

It is easily checked that any proper short exact sequence
(\ref{3}) of $G\rtimes \Gamma$-modules induces long exact homology
and cohomology sequences
\begin{multline*}
\cdots\rightarrow H^{\Gamma}_{n+1}(G,C_{2})\rightarrow
H^{\Gamma}_{n}(G,C_{1})\rightarrow \cdots\rightarrow
H^{\Gamma}_{2}(G,C_{2})\rightarrow
H^{\Gamma}_{1}(G,C_{1})\rightarrow \\ H^{\Gamma}_{1}(G,C)
\rightarrow H^{\Gamma}_{1}(G,C_{2})\rightarrow
H^{\Gamma}_{0}(G,C_{1})\rightarrow H^{\Gamma}_{0}(G,C)\rightarrow
H^{\Gamma}_{0}(G,C_{2})\rightarrow 0\;,
\end{multline*}

\begin{multline*}
0\rightarrow H^{0}_{\Gamma}(G,C_{1})\rightarrow
H^{0}_{\Gamma}(G,C)\rightarrow H^{0}_{\Gamma}(G,C_{2})\rightarrow
H^{1}_{\Gamma}(G,C_{1})\rightarrow
H^{1}_{\Gamma}(G,C)\rightarrow\\
H^{1}_{\Gamma}(G,C_{2})\rightarrow
H^{2}_{\Gamma}(G,C_{1})\rightarrow\cdots\rightarrow
H^{n}_{\Gamma}(G,C_{2})\rightarrow
H^{n+1}_{\Gamma}(G,C_{1})\rightarrow\cdots\;.
\end{multline*}

\

\

\section{Equivariant (co)homology of groups as cotriple
(co)homology}

\

To present the equivariant (co)homology of groups as cotriple
(co)homology we will use the free cotriple defined in the category
$\mathcal{G}_{\Gamma}$ of $\Gamma$-groups given in \cite{InN,
InHInN} to develop a non-abelian homology theory of groups. This
cotriple corresponds to the tripleability of
$\mathcal{G}_{\Gamma}$ over $\Gamma$-Sets. The resulting cotriple
$\mathcal{F}=(F,\tau,\delta)$ is the free cotriple in the category
of groups endowed with the $\Gamma$-action defined as follows. For
any $\Gamma$-group $G$ the action of $\Gamma$ on the free group
$F(G)$ is given by $^{\sigma}|g|= |^{\sigma}g|$, $g\in G$,
$\sigma\in \Gamma$. The cotriple thus defined essentially differs
from the cotriple introduced in \cite{CGO2} for the cotriple
interpretation of the cohomology of groups with operators. Let
$\mathcal{P}_{\mathcal{F}}$ be the projective class induced by the
cotriple $\mathcal{F}$ in the category $\mathcal{G}_{\Gamma}$. It
is easy to see that a morphism $f: G\rightarrow H$ of
$\Gamma$-groups is a $\mathcal{P}_{\mathcal{F}}$-epimorphism if it
is surjective and $\Gamma$-splitting. Since the category
$\mathcal{G}_{\Gamma}$ has finite limits, any $\Gamma$-group $G$
has a $\mathcal{P}_{\mathcal{F}}$-projective resolution
$(X_{*},\partial^{0}_{0},G)$ in the category
$\mathcal{G}_{\Gamma}$ in the sense of \cite{TiVo}, that is
$X_{*}$ is an augmented  pseudo-simplicial $\Gamma$-group which is
$\mathcal{P}_{\mathcal{F}}$-exact \cite{TiVo, InH2} and each
$X_{n}$, $n\geq 0$, belongs to the class
$\mathcal{P}_{\mathcal{F}}$. Many examples of pseudo-simplicial
sets which are not simplicial are given in \cite{Lo1, Fr1, Fr2}. A
$\mathcal{P}_{\mathcal{F}}$-epimorphism $f: P\rightarrow G$ with P
an object of the class $\mathcal{P}_{\mathcal{F}}$ will be called
a projective presentation of the $\Gamma$-group $G$. Any
projective presentation of $G$ induces in a natural way a
$\mathcal{P}_{\mathcal{F}}$-projective resolution
$P_{*}\rightarrow G$, constructed as follows:
$$
\cdots \;\dtwo\; F(L_{2}G)\overset{\tau_{L_{2}G}}{\rightarrow}
L_{2}G\underset{l^{2}_{2}}{\overset{l^{2}_{0}}{\;\three\;}}
F(L_{1}G)\overset{\tau_{L_{1}G}}{\rightarrow} L_{1}G
\underset{l^{1}_{1}}{\overset{l^{1}_{0}}{\;\two \;}}
P\overset{f}{\to} G ,
$$
where $(L_{1}G,l^{1}_{0},l^{1}_{1})$ is the simplicial kernel of
the morphism $f$, $(L_{2}G,l^{2}_{0},l^{2}_{1},l^{2}_{2})$ the
simplicial kernel of the pair of morphisms
$(l^{1}_{0}\tau_{L_{1}G},l^{1}_{1}\tau_{L_{1}G})$ and if
$(L_{n}G,l^{n}_{0},\ldots,l^{n}_{n})$ has been constructed, then
$(L_{n+1}G,l^{n+1}_{0},\ldots,l^{n+1}_{n+1})$ is the simplicial
kernel of the sequence of morphisms
$(l^{n}_{0}\tau_{L_{n}G},\ldots,l^{n}_{n}\tau_{L_{n}G})$.
Simplicial kernels are defined in \cite{TiVo, InH2}. Denote
$P_{0}= P$, $P_{n}= F(L_{n}(G))$ and $\partial_{i}^{n}=
l_{i}^{n}\tau_{L_{n}(G)}$ for $n> 0$.

Let $T$ be a functor from the category $\mathcal{G}_{\Gamma}$ to
the category $\mathcal{G}$ of groups.  Then by definition the left
cotriple derived functors $L^{\mathcal{F}}_{n}T(G)$ are equal to
$\pi_{n}(TF_{*}(G)$, $n\geq 0$, where $\tau: F_{*}(G)\rightarrow
G$ is the free cotriple resolution of $G$ with $F_{0}= F(G),
F_{n}= F(F^{n}(G)), n\geq 1$, $\partial_{i}^{n}= F^{i}\tau
F^{n-i}, s_{i}^{n}= F^{i}\delta F^{n-i}$(see \cite{Sw2}) and the
left derived functors $L_{n}^{\mathcal{P}_{\mathcal{F}}}T(G)$ with
respect to the projective class $\mathcal{P}_{\mathcal{F}}$ are
equal to $\pi_{n}(T(X_{*}))$, $n\geq 0$ (see \cite{InH1, InH2}).
If the functor $T$ is a contravariant functor with values in the
category of abelian groups, one has also its right derived
functors.

\medskip

\begin{prop}
The left cotriple derived functors of a functor $T:
\mathcal{G}_{\Gamma}\rightarrow \mathcal{G}$ are isomorphic to its
left $\mathcal{P}_{\mathcal{F}}$-derived functors.
\end{prop}
\begin{proof}
Let $(P_{*},\partial^{0}_{0},G)$ be the standard
$\mathcal{P}_{\mathcal{F}}$-resolution of $G$. Then it is easy to
see that this resolution is left contractible in the category of
$\Gamma$-Sets. Thus the augmented pseudosimplicial $\Gamma$-groups
$(F_{i}(P_{*}),F_{i}(\partial^{0}_{0}),F_{i}(G))$ are left
contractible for $i\geq 0$ in the category of $\Gamma$-groups (for
the categorical definition of left contractibility see \cite{Sw2,
InH2}). On the other hand the augmented simplicial $\Gamma$-groups
$(F_{*}(P_{j}),\tau_{j},P_{j})$ are also left contractible for
$j\geq 0$. It follows that the homotopy groups $\pi_{n}$, $n>0$,
of the pseudosimplicial groups $TF_{i}(P_{*})$ and $TF_{*}(P_{j})$
for $i,j\geq 0$ are trivial and the homotopy groups of
$\pi_{0}(TF_{i}(P_{*}))$ and $\pi_{0}(TF_{*}(P_{j}))$ give the
left projective and cotriple derived functors respectively of the
functor $T$. Consider now the bipseudosimplicial group $G_{**}(G)$
by putting $G_{pq}(G)= TF_{p}(P_{q}(G))$ and apply the Quillen
spectral sequences \cite{Qu, InH1, InH2} for a bipseudosimplicial
group. It follows that the $n$-th homotopy groups, $n\geq 0$, of
$TP_{*}(G)$ and $TF_{*}(G)$ are both isomorphic to the $n$-th
homotopy group of the diagonal pseudosimplicial group $\Delta
G_{**}$. It remains to apply Theorems 1.2 and 2.1 of \cite{InH1}
showing that the definition of the left projective derived
functors are independent of the projective resolution of $G$.
\end{proof}

\medskip

Note that this proposition is known for the left derived functors
of functors (right derived functors of contravariant functors)
with values in the category of abelian groups \cite{TiVo}.

Let $A$ be a fixed $\Gamma$-module and $_{A}\mathcal{G}_{\Gamma}$
the category of $\Gamma$-groups acting on $A$ such that the
condition (\ref{1}) holds. Consider the following functors from
the category $_{A}\mathcal{G}_{\Gamma}$ to the category of abelian
groups: $I(-)\otimes_{G\rtimes \Gamma} A$ and
$\Der_{\Gamma}(-,A)$, where $I(G)$ is the kernel of the canonical
homomorphism $\epsilon:Z(G)\rightarrow Z$ of $G\rtimes
\Gamma$-modules and $\Der_{\Gamma}(G,A)$ is the group of
$\Gamma$-derivations from $G$ to $A$ consisting of derivations
$f:G\rightarrow A$ such that $f(^{\sigma}g)= \; ^{\sigma}f(g)$,
$g\in G$, $\sigma\in \Gamma$ \cite{CGO2}.

\medskip

\begin{thm}\label{t-3}
There are isomorphisms
$$
H^{\Gamma}_{n}(G,A) \cong
L^{\mathcal{F}}_{n-1}(I(G)\otimes_{G\rtimes \Gamma} A)\;,\quad
H^{n}_{\Gamma}(G,A) \cong
R^{n-1}_{\mathcal{F}}\Der_{\Gamma}(G,A)\;,\quad n\geq 2.
$$
\end{thm}
\begin{proof} Apply the functor $I(-)$ to the free
cotriple resolution $\tau: F_{*}\rightarrow G$ of the
$\Gamma$-group $G$. One gets an augmented simplicial $G\rtimes
\Gamma$-module $I(F_{*})\rightarrow I(G)$. We introduce the
notations:
$$
IF_{n}(G)\cong \underset{y\in F^{n}(G)}{\sum}
\mathbb{Z}(F^{n}(G))(y-e)= D_{n}(G)\;,\quad \underset{y\in
F^{n}(G)}{\sum} \mathbb{Z}(G)(y-e)= E_{n}(G)
$$
for $n\geq 0$, $F^{0}(G) = G$.

There are natural homomorphisms
$$
\alpha_{n}: D_{n}(G)\rightarrow
E_{n}(G)\;, \quad n\geq 1\;,
$$
induced by the homomorphism
$\tau\partial^{1}_{0}\partial^{2}_{0}\cdots
\partial^{n-2}_{0}\partial^{n-1}_{0}: F^{n}(G)\rightarrow G$ such
that we obtain a morphism of augmented simplicial $G\rtimes
\Gamma$-modules
$$
(D_{*}(G)\rightarrow I(G))\rightarrow (E_{*}(G)\rightarrow
I(G))\;.
$$
The left $\Gamma$-contractibility of
the cotriple resolution $F_{*}(G)\rightarrow G$ implies the
$\Gamma$-contractibility of the corresponding induced abelian
chain complexes
$$
\cdots \rightarrow D_{n}(G)\rightarrow \cdots \rightarrow
D_{2}(G)\rightarrow D_{1}(G)\rightarrow I(G)\rightarrow 0,
$$
$$
\cdots \rightarrow E_{n}(G)\overset{\varepsilon_{n}}{\rightarrow}
\cdots \rightarrow E_{2}(G)\overset{\varepsilon_{2}}{\rightarrow}
E_{1}(G)\overset{\varepsilon_{1}}{\rightarrow} I(G)\rightarrow 0,
$$
where $\varepsilon_{n}=\underset{i}{\sum}
(-1)^{i}\varepsilon^{n}_{i}$, $n\geq 1$. In effect the canonical
$\Gamma$-injections $\{f, f_{n}, n\geq 1\}$,
$$
f: G\rightarrow F(G)\;, \quad f_{n}: F^{n}(G)\rightarrow
F^{n+1}(G)= (F(F^{n}(G))) \;,\quad n\geq 1\;,
$$
 yield the left
$\Gamma$-contractibility of $F_{*}(G)\rightarrow G$ in the
category of $\Gamma$-sets (see [33], Lemma 1.2). Therefore we
obtain $\Gamma$-homomorphisms
$$
\mathbb{Z}(f): \mathbb{Z}(G)\rightarrow \mathbb{Z}(F(G))\;,\quad
\mathbb{\mathbb{Z}}(f_{n}): \mathbb{Z}(F^{n}(G))\rightarrow
\mathbb{Z}(F^{n+1}(G))\;,\; n\geq 1\;,
$$
of free abelian $\Gamma$- groups induced by $\{f, f_{n}, n\geq
1\}$, where the action of $\Gamma$ on $\mathbb{Z}(F^{n}(G))$,
$n\geq 0$, $F^{0}(G)= G$, is induced by the above defined action
of $\Gamma$ on $F^{n}(G)$. The action of $\Gamma$ on
$IF_{n}(G)\cong \underset{y\in F^n(G)}{\Sigma}
\mathbb{Z}(F^{n}(G))(y-e), n\geq 0$, is induced by the action of
$\Gamma$ on $\mathbb{Z}(F^{n}(G))$, namely
$$
^{\sigma}(x(y-e))\;= \; ^{\sigma}x (^{\sigma}y-e)\;,\quad x\in
\mathbb{Z}(F^{n}(G))\;,\; y\in F^{n}(G)\;.
$$
The $\Gamma$- homomorphisms $\{\mathbb{Z}(f), \mathbb{Z}(f_{n}),
n\geq 0\}$ induce $\Gamma$- homomorphisms
$$
IG\rightarrow IF_{0}(G)\;,\quad IF^{n}(G)\rightarrow
IF^{n+1}(G)\;,\; n\geq 0\;,
$$
and yield the required $\Gamma$- contraction in $IF_{*}(G)$.

Thus each short exact sequence

$$
0\rightarrow \Ker\varepsilon_{n}\rightarrow E_{n}(G)\rightarrow
Im\varepsilon_{n}\rightarrow 0,\;\; n\geq 1,
$$
is $\Gamma$-splitting and it follows that $(E_{*}(G)\rightarrow
I(G))$ is a relatively free resolution of the $G\rtimes
\Gamma$-module $I(G)$.

It is obvious that the homomorphisms $\alpha_{n}$, $n\geq 1$,
induce isomorphisms
\begin{equation}
\begin{matrix}\label{6}
D_{n}(G)\otimes_{F^{n}(G)\rtimes \Gamma} A \cong
E_{n}(G)\otimes_{G\rtimes \Gamma} A\;,\\
\Hom_{F^{n}(G)\rtimes\Gamma}(D_{n}(G),A) \cong
\Hom_{G\rtimes\Gamma}(E_{n}(G),A)\;.
\end{matrix}
\end{equation}
Whence we deduce from (\ref{6}) that
$L^{\mathcal{F}}_{n}(I(G)\otimes A) \cong
\Tor^{\mathcal{P}}_{n}(I(G),A)$, $n\geq 0$. It is easily checked
that the well known isomorphism
$$
\Hom_{F^{n}(G)}(D_{n}(G),A)\cong \Der(F^{n}(G),A)
$$
is compatible with the action of $\Gamma$, whence its restriction
on the subgroup \linebreak $\Hom_{F^{n}(G)\rtimes
\Gamma}(D_{n}(G),A)$ gives an isomorphism with
$\Der_{\Gamma}(F^{n}(G),A)$. Thus from (\ref{6}) one gets
$R^{n}_{\mathcal{F}}\Der_{\Gamma}(G,A) \cong
\Ext^{n}_{\mathcal{P}}(I(G),A)$, $n\geq 0$.

The proper short exact sequence of $G{\rtimes\Gamma}$-modules
\begin{align}\label{7}
0\rightarrow I(G)\rightarrow Z(G)\rightarrow Z\rightarrow 0
\end{align}
yields long exact sequences of the relative derived functors of
the functors $-\otimes_{G\rtimes \Gamma} A$ and $\Hom_{G\rtimes
\Gamma}(-,A)$ implying the isomorphisms
$$
H^{\Gamma}_{n+1}(G,A)\cong \Tor^{\mathcal{P}}_{n}(I(G),A) \quad
\text{and} \quad H^{n+1}_{\Gamma}(G,A) \cong
\Ext^{n}_{\mathcal{P}}(I(G),A)\;,\; n\geq 1\;,
$$ which give the required isomorphisms.
\end{proof}

\medskip

It is clear that $L^{\mathcal{F}}_{0}(I(G)\otimes_{G\rtimes
\Gamma} A) \cong I(G)\otimes_{G\rtimes \Gamma} A$ and
$R^{0}_{\mathcal{F}}\Der_{\Gamma}(G,A) \cong \Der_{\Gamma}(G,A)$.

\medskip

\begin{defn}A $\Gamma$-group $G$ will be called $\Gamma$-free, if it
is a free group with basis a $\Gamma$-subset.
\end{defn}

\medskip

\begin{cor}If $G$ is a retract of a $\Gamma$-free group, then $H^{\Gamma}_{n}(G,A)=
0$ and $H^{n}_{\Gamma}(G,A) = 0$ for $n> 1$ and any $G\rtimes
\Gamma$-module $A$.
\end{cor}
\begin{proof} The augmented simplicial group
$F_{*}(G)\rightarrow G$ is left contractible \cite{Sw2} implying
the triviality of the homotopy groups
$\pi_{n}(IF_{*}(G)\otimes_{G\rtimes \Gamma} A)$ and
$\pi_{n}\Der_{\Gamma}(F_{*}(G),A)$ for $n\geq 1$.
\end{proof}

\medskip

\begin{cor}\label{c-6}
If $G$ and $\Gamma$ act trivially on $\mathbb{Z}$, then
$H^{\Gamma}_{1}(G,\mathbb{Z}) \cong I(G)\otimes_{G\rtimes\Gamma}
\mathbb{Z}$.
\end{cor}
\begin{proof} The proof follows immediately from the long exact
sequence of the functors $\Tor^{\mathcal{P}}_{n}(-,\mathbb{Z})$
induced by (\ref{7}), since $\mathbb{Z}(G)\otimes_{G\rtimes
\Gamma} \mathbb{Z} \cong \mathbb{Z}$.
\end{proof}

\medskip

\begin{prop}\label{p-7}
The cotriple derived functors
$L^{\mathcal{F}}_{n}H^{\Gamma}_{1}(-,A)$ are isomorphic to
\linebreak $H^{\Gamma}_{n+1}(-,A)$, $n> 0$.
\end{prop}
\begin{proof} The long exact sequence of the functors
$\Tor^{\mathcal{F}}_{n}(-,A)$ for the sequence (\ref{7}) yields
the exact sequence
$$
0\rightarrow H^{\Gamma}_{1}(G,A)\rightarrow I(G)\otimes_{G\rtimes
\Gamma} A\rightarrow \mathbb{Z}(G)\otimes_{G\rtimes \Gamma}
A\rightarrow \mathbb{Z}\otimes_{G\rtimes \Gamma} A\rightarrow 0.
$$
It follows that there is a functorial short exact sequence
$$
0\rightarrow H^{\Gamma}_{1}(G,A)\rightarrow I(G)\otimes_{G\rtimes
\Gamma } A\rightarrow I_{G\rtimes \Gamma}A /
I_{\Gamma}A\rightarrow 0
$$
inducing a short exact sequence of abelian simplicial groups
$$
0\rightarrow H^{\Gamma}_{1}(F_{*}(G),A)\rightarrow
I(F_{*}(G))\otimes_{G\rtimes \Gamma} A\rightarrow
I_{F_{*}(G)\rtimes \Gamma} A / I_{\Gamma}A\rightarrow 0.
$$
It remains to apply the corresponding long exact homotopy sequence
and to see that $I_{F_{*}(G)\rtimes \Gamma}A / I_{\Gamma}A$ is a
constant abelian simplicial group.
\end{proof}

\medskip

Denote by $[G,G]_{\Gamma}$ the subgroup of the $\Gamma$-group $G$
generated  by $[G,G]$ and by the elements of the form
$^{\sigma}g\cdot g^{-1}$, $g\in G$, $\sigma\in \Gamma$. This
subgroup will be called the $\Gamma$-commutator subgroup of $G$.
It is obvious that $[G,G]_{\Gamma}$ is a normal $\Gamma$-subgroup
of $G$ and $\Gamma$ acts trivially on the abelian group $G /
[G,G]_{\Gamma}$. If $H$ is a normal $\Gamma$-subgroup of $G$, we
denote by $[G,H]_{\Gamma}$ the subgroup of $G$ generated by the
elements $x\;^{\sigma}yx^{-1}y^{-1}$, where $x\in G, y\in H,
\sigma\in \Gamma$.

Let $B$ be an abelian group on which $\Gamma$ acts trivially and
$f:G\rightarrow B$ a homomorphism of $\Gamma$-groups. Then $f$
factorizes uniquely through $G / [G,G]_{\Gamma}$. Consider the
subgroup of $G$ generated by the elements of the form
$x\;^{\sigma}yx^{-1}y^{-1}$, $x,y\in G$, $\sigma\in \Gamma$. It is
easily seen that this subgroup coincides with $[G,G]_{\Gamma}$.
The elements $x\;^{\sigma}yx^{-1}y^{-1}= [x,y]_{\sigma}$ will be
called $\Gamma$-commutators of $G$, the group $G / [G,G]_{\Gamma}$
the $\Gamma$-abelianization $G^{ab}_{\Gamma}$ of $G$ and the
corresponding functor $G\mapsto G^{ab}_{\Gamma}$ the
$\Gamma$-abelianization functor.

\medskip

\begin{prop}\label{p-8}
There is a functorial isomorphism
$$
I(G)\otimes_{G\rtimes \Gamma}A \cong G / [G,G]_{\Gamma}\otimes
A\;,
$$
where $G$ and $\Gamma$ act trivially on $A$.
\end{prop}
\begin{proof} It is enough to show that
$I(G)\otimes_{G\rtimes \Gamma} Z \cong G / [G,G]_{\Gamma}$. This
isomorphism is given by $(g-e)\otimes n \mapsto n[g]$. Its
converse is defined by $[g]\mapsto (g-e)\otimes 1$. We have only
to show the correctness of the converse map.

One has
\begin{align*}
&(x\;^{\sigma}yx^{-1}y^{-1}-e)\otimes 1 =
(x\;^{\sigma}yx^{-1}y^{-1}-x+x-e)\otimes 1 =
x(^{\sigma}yx^{-1}y^{-1}-e)\otimes 1 + \\
&(x-e)\otimes 1 = (^{\sigma}yx^{-1}y^{-1}-e)\otimes 1 +
(x-e)\otimes 1 = (^{\sigma}yx^{-1}y^{-1}-^{\sigma}y +
^{\sigma}y-e)\otimes 1 +\\
&(x-e)\otimes 1 = \;^{\sigma}y(x^{-1}y^{-1}-e)\otimes 1 +
(^{\sigma}y-e)\otimes 1 + (x-e)\otimes 1 = (x^{-1}y^{-1}-e)\otimes 1 + \\
&(y-e)\otimes 1 + (x-e)\otimes 1 =
(x^{-1}y^{-1}-x^{-1}+x^{-1}-e)\otimes 1 + (y-e)\otimes 1 +
(x-e)\otimes 1 =\\& x^{-1}(y^{-1}-e)\otimes 1 + (x^{-1}-e)\otimes
1 + (y-e)\otimes 1 + (x-e)\otimes 1 = (y^{-1}-e)\otimes 1 +\\&
(x^{-1}-e)\otimes 1 + (y-e)\otimes 1 + (x-e)\otimes 1.
\end{align*}
But for any element $x\in G$ the following equalities hold:
$$0=0\otimes 1 = (xx^{-1}-e)\otimes 1= (x-e)\otimes 1 +
(x^{-1}-e)\otimes 1\;.$$ It follows that under the afore defined
converse map any $\Gamma$-commutator becomes 0 showing its
correctness.
\end{proof}

The isomorphism of Proposition 8 holds only for trivial actions on
$A$ and it is natural in the sense that in this case it is
functorial and in fact uniquely defined. For an arbitrary
$G$-module in the classical case there exists another form of this
isomorphism, where its right side is replaced by the non-abelian
tensor product of the groups $G$ and $A$.

\

\

\section{The equivariant integral homology}

\

Denote the $\Gamma$-equivariant integral homology groups
$H^{\Gamma}_{n}(G,\mathbb{Z}) = H^{\Gamma}_{n}(G)$, $n\geq 0$, the
groups $G$ and $\Gamma$ acting trivially on $\mathbb{Z}$.

\medskip

\begin{cor}\label{c-9}
There is a functorial isomorphism
$$
H^{\Gamma}_{1}(G) \cong G /[G,G]_{\Gamma}\;.
$$
\end{cor}
\begin{proof} The proof follows from Corollary \ref{c-6} and Proposition \ref{p-8}.
\end{proof}

\medskip

Note that the group $[G,G]_{\Gamma} / [G,G]$ shows the difference
between the classical and the $\Gamma$-equivariant abelianization
functors. We denote by $T\Gamma$ the functor assigning to any
$\Gamma$-group $G$ the abelian group $[G,G]_{\Gamma} / [G,G]$. We
also denote by $\Gamma\cdot G$ the subgroup of $G$ generated by
the elements $^{\sigma}g\cdot g^{-1}$, $g\in G$, $\sigma\in
\Gamma$. One has a natural homomorphism
$$\beta_{n}:H_{n}(G)\rightarrow H^{\Gamma}_{n}(G)\;,\quad n\geq 0\;,
$$
induced by the morphism of abelian simplicial groups
$$(I(F_{*}(G))\otimes_{ G} \mathbb{Z})\rightarrow
I(F_{*}(G))\otimes_{G\rtimes \Gamma} \mathbb{Z}\;.$$

\medskip

\begin{thm}\label{t-10}
\begin{enumerate}
\item[(i)] There is an isomorphism
$$
L^{\mathcal{F}}_{n}(G^{ab}_{\Gamma})\cong
H^{\Gamma}_{n+1}(G)\;,\quad n\geq 0\;.
$$
\item[(ii)] There are a
functorial short exact sequence
$$
0\rightarrow \Gamma\cdot G / [G,G]\cap \Gamma\cdot G\rightarrow
H_{1}(G)\rightarrow H^{\Gamma}_{1}(G)\rightarrow 0\;,
$$
and a long exact homology sequence
\begin{multline*}
\cdots \rightarrow H^{\Gamma}_{n+1}(G)\rightarrow
L^{\mathcal{F}}_{n-1}T\Gamma(G)\rightarrow H_{n}(G)\rightarrow
H^{\Gamma}_{n}(G)\rightarrow
L^{\mathcal{F}}_{n-2}T\Gamma(G)\rightarrow \cdots \rightarrow\\
H^{\Gamma}_{3}(G)\rightarrow
L^{\mathcal{F}}_{1}T\Gamma(G)\rightarrow H_{2}(G)\rightarrow
H^{\Gamma}_{2}(G)\rightarrow
L^{\mathcal{F}}_{0}T\Gamma(G)\rightarrow H_{1}(G)\rightarrow
H^{\Gamma}_{1}(G)\rightarrow 0.
\end{multline*}
\end{enumerate}
\end{thm}
\begin{proof}(i) The proof follows from Proposition \ref{p-7} and Corollary \ref{c-9}.

(ii) The commutative diagram
$$
\begin{matrix}
H_{1}(G)& \overset{\beta_1}{\rightarrow} & H_{1}^{\Gamma}(G)\\
\vcenter{\llap{$_{\cong}{}$}}\downarrow & & \downarrow\vcenter{\rlap{${}_{\cong}$}} \\
G/[G,G] & \rightarrow & G/[G,G]_\Gamma
\end{matrix}
$$
shows that $\Ker\beta_{1}$ is isomorphic to $[G,G]_{\Gamma} /
[G,G]$ and it is clear that $[G,G]_{\Gamma}= [G,G]\cdot
(T\Gamma(G))$. By applying the short exact sequence (ii) to the
cotriple resolution $F_{*}(G)\rightarrow G$, we obtain a short
exact sequence of simplicial abelian groups
$$
0\rightarrow T\Gamma(F_{*}(G))\rightarrow
H_{1}(F_{*}(G))\rightarrow H^{\Gamma}_{1}(F_{*}(G))\rightarrow 0
$$
inducing a long exact homology sequence and it remains to recall
that the cotriple $n$-th derived functor of the first integral
homology gives the $(n+1)$-th integral homology group, $n\geq 0$.
\end{proof}

\medskip

Let $A$ be a $G\rtimes \Gamma$-module on which $G$ and $\Gamma$
act trivially. Then
$$
H^{\Gamma}_{1}(G,A)\cong I(G)\otimes_{G\rtimes\Gamma} A \cong
H^{\Gamma}_{1}(G)\otimes A\cong G / [G,G]_{\Gamma}\otimes A
$$
and
$$
H^{1}_{\Gamma}(G,A)\cong \Der_{\Gamma}(G,A)\cong \Hom_{G\rtimes
\Gamma}(I(G)\otimes_{G\rtimes \Gamma} \mathbb{Z},A)\cong \Hom(G /
[G,G]_{\Gamma},A)\;.
$$
On the other hand, if $G$ is a $\Gamma$-free
group with basis $X$ and $\beta: G\rightarrow G / [G,G]_{\Gamma}$
is the canonical $\Gamma$-homomorphism, then $G / [G,G]_{\Gamma}$
is a free abelian group with basis $\beta(X)$. Any map
$\gamma:\beta(X)\rightarrow B$ to an abelian group $B$ induces a
$\Gamma$-map $\gamma\beta:X\rightarrow B$ which is uniquely
extended to a $\Gamma$-homomorphism $G\rightarrow B$ assuming
$\Gamma$ acts trivially on $B$ and one gets a uniquely defined
homomorphism $G / [G,G]_{\Gamma}\rightarrow B$ whose restriction
on $\beta(X)$ is equal to $\gamma$.

We deduce that for any $G\rtimes \Gamma$-module $A$ with trivial
actions of $G$ and $\Gamma$ on $A$ we obtain universal coefficient
formulas for the equivariant (co)homology groups
$H^{\Gamma}_{n}(G,A)$ and $H^{n}_{\Gamma}(G,A)$, $n\geq 0$.

\medskip

\begin{thm}
There are short exact split (not naturally) sequences
\begin{align*}
&0\rightarrow H^{\Gamma}_{n} (G)\otimes A\rightarrow
H^{\Gamma}_{n}(G,A)\rightarrow
\Tor_{1}(H^{\Gamma}_{n-1}(G),A))\rightarrow 0\;,\\ & 0\rightarrow
\Ext^{1}(H^{\Gamma}_{n-1}(G),A))\rightarrow
H^{n}_{\Gamma}(G,A)\rightarrow
\Hom(H^{\Gamma}_{n}(G),A))\rightarrow 0
\end{align*}
for $n\geq 0$.
\end{thm}

\

\

\section{Universal central $\Gamma$-equivariant extensions and Hopf
formula}

\

Let $G$ be a $\Gamma$-group and $A$ a $G\rtimes \Gamma$-module.

\medskip

\begin{defn}
A $\Gamma$-equivariant extension $E$ of the $\Gamma$-group $G$ by
the $\Gamma$-equivariant $G$-module $A$ is an extension of $G$ by
the $G$-module $A$
$$
E: 0\rightarrow A\overset{\alpha}{\rightarrow}
B\overset{\beta}{\rightarrow} G\rightarrow 1
$$
satisfying the following conditions:
\begin{enumerate}
\item[1)] $E$ is a sequence of $\Gamma$-groups,
\item[2)] $E$ is $\Gamma$-splitting, that is there is a $\Gamma$-map
$\gamma:G\rightarrow  B$ such that $\beta\gamma = 1_{G}$.
\end{enumerate}
$E$ is called a central $\Gamma$-equivariant extension of the
$\Gamma$-group $G$, if $\alpha(A)$ belongs to the center of $B$
and $\Gamma$ acts trivially on $A$.
\end{defn}

\medskip

Note that $\Gamma$-equivariant extensions of $\Gamma$-groups
investigated in \cite{CGO2} do not in general satisfy the
condition 2).

A central $\Gamma$-equivariant extension $(U,\beta)$ of $G$ is
called universal, if for any central $\Gamma$-equivariant
extension $(X,\alpha)$ of $G$ there is a unique
$\Gamma$-homomorphism $U\rightarrow X$ over $G$.

Two $\Gamma$-equivariant extensions $E$ and $E'$ of $G$ by $A$ are
called equivalent if there is a morphism $E\rightarrow E'$ which
is the identity on $A$ and $G$. We denote by $E_{\Gamma}(G,A)$ the
set of equivalence classes of $\Gamma$-equivariant extensions of
$G$ by $A$.

\medskip

\begin{defn} A $\Gamma$-group $G$ is called $\Gamma$-perfect, if
$G$ coincides with its $\Gamma$-commutator subgroup
$[G,G]_{\Gamma}$ (see also \cite{Lo2}).
\end{defn}

\medskip

Below we give important examples of $\Gamma$-groups which are
$\Gamma$-perfect but not perfect (see Section 6).

\medskip

\begin{prop}\label{p-14}
If $(X,\varphi)$ is a central $\Gamma$-equivariant extension of a
$\Gamma$-perfect group $G$, then the $\Gamma$-commutator subgroup
$X'= [G,G]_{\Gamma}$ is $\Gamma$-perfect and maps onto $G$.
\end{prop}
\begin{proof} Since $G$ is $\Gamma$-perfect, it is
clear that $\varphi$ maps $X'$ onto $G$. It follows that any
element $x\in X$ can be written as a product $x'c$ with $x'\in X'$
and $c$ belongs to $Ker\varphi$. Therefore every generator of $X'$
of the form $[x_{1},x_{2}]$ is equal to
$[x'_{1}c'_{1},x'_{2}c'_{2}]\;=\; [x'_{1},x'_{2}]$ with $
x'_{1},x'_{2}\in X'$ and of the form $^{\sigma}x\cdot x^{-1}$ is
equal to ${}^{\sigma}(x'c)\cdot (x'c)^{-1}\;={}^{\sigma}x'\cdot
x'^{-1}$ with $x'\in X'$. Whence $X'= [X',X']_{\Gamma}$.
\end{proof}

\medskip

\begin{cor}
$(X',\varphi|_{X'})$ is a central $\Gamma$-equivariant
extension of $G$.
\end{cor}

\begin{proof} We have only to show that
$(X',\varphi|_{X'})$ is $\Gamma$-splitting. Let $\gamma$ be the
$\Gamma$-splitting map for $(X,\varphi)$, that is $\varphi\gamma=
1_{G}$. Then a $\Gamma$-splitting map $\gamma'$ for the extension
$(X',\varphi|_{X'})$ is defined as follows: consider the
decomposition $\{D_{\eta}\}$ of $G$ into orbits with respect to
the action of $\Gamma$ on $G$. Choose a representative
$z_{\eta}\in D_{\eta}$ for each $D_{\eta}$ and choose an
expression
$$
z_{\eta}= [x_{1},y_{1}]_{{\sigma}_{1}}\cdot
[x_{2},y_{2}]_{{\sigma}_{2}}\cdots [x_{k},y_{k}]_{{\sigma}_{k}}
$$
of $z_{\eta}$ in terms of $\Gamma$- commutators. The equalities
\begin{align*}
&{}^{{\sigma}'}{([x,y]_{\sigma})}\;={}^{{\sigma}'}{( x\;^{\sigma}
y x^{-1} y^{-1})}=\; {{^{\sigma}}'x\;{^{\sigma}}'^{\sigma}}y
\;{^{\sigma}}'
(x^{-1}) \;{^{\sigma}}'(y^{-1})=   \\
&{}{^{\sigma}}'x\;
{^{\sigma}}'^{\sigma}{^{\sigma}}'^{-1}({^{\sigma}}'y)
\;{^{\sigma}}'(x^{-1})\;{^{\sigma}}'(y^{-1})=
[{^{\sigma}}'x,\;{^{\sigma}}'y]_{\sigma ' \sigma {\sigma
'}^{-1}}\;,
\end{align*}
$x,y\in G$, $\sigma,\sigma'\in \Gamma$, imply the expression
$$
^{\sigma}z_{\eta}= [^{\sigma}x_1,\;^{\sigma}y_1]_{\sigma \sigma_1
\sigma^{-1}}\cdot [^{\sigma}x_2,\;^{\sigma}y_2]_{\sigma \sigma_2
\sigma^{-1}}\cdots [^{\sigma}x_k,\;^{\sigma}y_k]_{\sigma
\sigma_{k} \sigma^{-1}}
$$
for $^{\sigma}z_{\eta}$, $\sigma\in \Gamma$.

Clearly any element $g\in G$ has the form $^{\sigma}z_{\eta}$ for
some $z_{\eta}\in D_{\eta}$ and $\sigma\in \Gamma$. The required
$\Gamma$- splitting map $\gamma': G\rightarrow X'$ is given by
setting
$$
\gamma'(g)= \gamma'(^{\sigma}z_{\eta})=
[^{\sigma}\gamma(x_1),\;^{\sigma}\gamma(y_1)]_{\sigma \sigma_1
\sigma^{-1}}\cdot
[^{\sigma}\gamma(x_2),\;^{\sigma}\gamma(y_2)]_{\sigma \sigma_2
\sigma^{-1}}\cdots
[^{\sigma}\gamma(x_k),\;^{\sigma}\gamma(y_k)]_{\sigma \sigma_{k}
\sigma ^{-1}}.
$$
\end{proof}

\medskip

\begin{thm}\label{t-16}
A central $\Gamma$-equivariant extension $(U,\beta)$ of a
$\Gamma$-group $G$ is universal if and only if $U$ is
$\Gamma$-perfect  and every central equivariant $\Gamma$-extension
$(W,\alpha)$ of $U$ splits.
\end{thm}
\begin{proof} We will follow Milnor's proof of the classical
case \cite{Mi}. Let $U$ be $\Gamma$-perfect and every central
$\Gamma$-equivariant extension of $U$ splits. Let $(X,\varphi)$ be
an arbitrary central $\Gamma$-equivariant extension of $G$. Take
the following diagram with exact rows
$$
\begin{matrix}
0&\rightarrow & C& \rightarrow & X \times U
&\overset{p}{\rightarrow}& U & \rightarrow & 1\\
& & || & & \vcenter{\llap{$_{q}{}$}}\downarrow & &
\downarrow\vcenter{\rlap{${}_{\beta}$}} & & \\
0&\rightarrow & C& \rightarrow & X
&\overset{\varphi}{\rightarrow}& G & \rightarrow & 1 ,
\end{matrix}
$$
where $X\times U$ is the fiber product of $X\rightarrow
G\leftarrow U$. It is easy to check that the top row is a central
$\Gamma$-equivariant extension of $U$. Therefore one has a
$\Gamma$-section $s:U\rightarrow X\times U$. Thus the
$\Gamma$-homomorphism $f=qs: U\rightarrow X$ is over $G$ and the
diagram
$$
\begin{matrix}
0&\rightarrow & D& \rightarrow & U
&\overset{\beta}{\rightarrow}& G & \rightarrow & 1\\
& & \downarrow & & \vcenter{\llap{$_{f}{}$}}\downarrow & &
|| & & \\
0&\rightarrow & C& \rightarrow & X
&\overset{\varphi}{\rightarrow}& G & \rightarrow & 1
\end{matrix}
$$
is commutative. To prove that $(U,\beta)$ is universal, it remains
to show the uniqueness of such an $f$. Let
$f_{1},f_{2}:U\rightarrow X$ be two $\Gamma$-homomorphisms over
$G$. Then one gets a $\Gamma$-homomorphism $h:U\rightarrow C$
given by $h(u)= f_{1}(u)\cdot f_{2}(u)^{-1}$, $u\in U$, which is
trivial, since $U$ is $\Gamma$-perfect and $\Gamma$ acts trivially
on $C$.

Let $(X,\varphi)$ and $(Y,\psi)$ be central $\Gamma$-equivariant
extensions of $G$. Then, as we have seen, if $Y$ is
$\Gamma$-perfect there exists at most one $\Gamma$-homomorphism
from $Y$ to $X$ over $G$. If $Y$ is not $\Gamma$-perfect, then
there is a suitable central $\Gamma$-equivariant $(X,\varphi)$ of
$G$ such that there exists more than one $\Gamma$-homomorphism
from $Y$ to $X$ over $G$. Indeed in this case there exists a non
trivial $\Gamma$-homomorphism $f$ from $Y$ to some abelian group
$A$ on which $\Gamma$ acts trivially. Take the central
$\Gamma$-equivariant split extension
$$
0\rightarrow A\rightarrow A\times G\rightarrow G\rightarrow 1.
$$
Setting $f_{1}(y)= (0,\psi(y))$ and $f_{2}(y)=(f(y),\psi(y))$ one
gets two distinct $\Gamma$-homomor-phisms $f_{1}$ and $f_{2}$ from
$Y$ to $A\times G$ over $G$.

Now let $(U,\beta)$ be a universal central $\Gamma$-equivariant
extension of $G$ and $(W,\alpha)$ a central $\Gamma$-equivariant
extension of $U$. Since $(U,\beta)$ is a universal central
$\Gamma$-equivariant extension of $G$, the group $U$ is
$\Gamma$-perfect. We will show that $(W,\beta\alpha)$ is a central
$\Gamma$-equivariant extension of $G$. Take $x_{0}\in
\Ker\beta\alpha$. Then $\alpha (x_{0})$ belongs to the center of
$U$ and $\Gamma$ acts trivially on $\Ker\beta\alpha$. In effect,
if $\gamma:U\rightarrow W$ is a $\Gamma$-section of $(W,\alpha)$,
one has $^{\sigma}(x_{0}-\gamma\alpha(x_{0}))= x_{0}-\gamma\alpha
(x_{0})$; on the other hand, $^{\sigma}(x_{0}- \gamma\alpha
(x_{0}))$= $^{\sigma}x_{0}- \gamma\;^{\sigma}\alpha(x_{0}) =$
$^{\sigma}x_{0}- \gamma\alpha(x_{0})$. Whence $^{\sigma}x_{0}=
x_{0}$, $x_{0}\in \Ker\beta\alpha$, $\sigma\in \Gamma$. We obtain
a $\Gamma$-homomorphism $h:W\rightarrow W$ over $U$ given by
$h(x)= x_{0}xx_{0}^{-1}$, $x\in W$. Since $[W,W]_{\Gamma}$ is
$\Gamma$-perfect, the restriction of $h$ to the
$\Gamma$-commutator subgroup $[W,W]_{\Gamma}$ of $W$ is the
identity map. Thus $x_{0}$ commutes with elements of
$[W,W]_{\Gamma}$ implying $x_{0}$ belongs to the center of $W$,
since $W$ is generated by $[W,W]_{\Gamma}$ and $\Ker\alpha$.

We deduce that there is a unique morphism $(U,\beta)\rightarrow
(W,\beta\alpha)$ over $G$, since $(U,\beta)$ is universal.
Therefore the composite $\alpha k$ of the induced
$\Gamma$-homomorphism $k:U\rightarrow W$ over $G$ with $\alpha$ is
equal to the identity showing that $(W,\alpha)$ splits.
\end{proof}

\medskip

It should be noted that central $\Gamma$-equivariant extensions
without $\Gamma$-splitting property were used in \cite{Lo2} to
characterize universal central relative extensions of an
epimorphism $\nu :\Gamma\rightarrow Q$ of groups.

Let $\tau: P\rightarrow G$ be a projective presentation of the
$\Gamma$-group $G$ and $R$ denotes the kernel of $\tau$. Then the
$\Gamma$-homomorphism $\tau$ sends the normal $\Gamma$-subgroup
$[P,R]_{\Gamma}$ of $P$ to $1$ and therefore induces a
$\Gamma$-homomorphism $\tau':[P,P]_{\Gamma} /
[P,R]_{\Gamma}\rightarrow [G,G]_{\Gamma}$ which is surjective.

\medskip

\begin{thm}\label{t-17}
\begin{enumerate}
\item[(i)] If $G$ is $\Gamma$-perfect, then $([P,P]_{\Gamma} /
[P,R]_{\Gamma}, \tau')$ is a universal central
$\Gamma$-equivariant extension of $G$.
\item[(ii)] For any
$\Gamma$-group $G$ the group $(R\cap [P,P]_{\Gamma}) /
[P,R]_{\Gamma}$ is isomorphic to $H^{\Gamma}_{2}(G)$.
\end{enumerate}
\end{thm}
\begin{proof} (i)The extension $(P /
[P,R]_{\Gamma},\tau)$ is a central $\Gamma$-equivariant extension
of $G$. Clearly this extension is central. The group $\Gamma$ acts
trivially on $\Ker\tau= R / [P,R]_{\Gamma}$. Indeed the element
$^{\sigma}x\cdot x^{-1}$, $x\in R$, belongs to $[P,R]_{\Gamma}$
for any $\sigma\in \Gamma$.

By Proposition \ref{p-14} the group $[P,P]_{\Gamma} /
[P,R]_{\Gamma}$ is $\Gamma$-perfect and maps onto $G$. By
Corollary 15 it follows that the extension $[P,P]_{\Gamma} /
[P,R]_{\Gamma}\rightarrow G$ has a $\Gamma$- equivariant splitting
map. Let $(X,\psi)$ be a central $\Gamma$-equivariant extension of
$G$. Then there is a $\Gamma$-homomorphism $f:P\rightarrow X$ over
$G$. Since $(X,\psi)$ is a central $\Gamma$-equivariant extension
of $G$, it is easy to see that $f([P,R]_{\Gamma})= 1$. Therefore
the restriction of $f$ on $[P,P]_{\Gamma}$ induces a unique
$\Gamma$-homomorphism $[P,P]_{\Gamma} / [P,R]_{\Gamma}\rightarrow
X$ over $G$, since $[P,P]_{\Gamma} / [P,R]_{\Gamma}$ is a
$\Gamma$-perfect group. We deduce that the short exact sequence
$$
0\rightarrow (R\cap [P,P]_{\Gamma}) / [P,R]_{\Gamma}\rightarrow
[P,P]_{\Gamma} / [P,R]_{\Gamma}\rightarrow G\rightarrow 1
$$
is a universal central $\Gamma$-extension of the $\Gamma$-perfect
group $G$.

(ii) Let $P_{*}\rightarrow G$ be the
$\mathcal{P}_{\mathcal{F}}$-projective resolution of the $\Gamma$-
group $G$ induced by the projective presentation
$\tau:P\rightarrow G$, which we have defined in section 2. The
long exact homotopy sequence induced by the short exact sequence
$$
1\rightarrow ([P_{*},P_{*}]_{\Gamma}\rightarrow
[G,G]_{\Gamma})\rightarrow (P_{*}\rightarrow G)\rightarrow
((P_{*})^{ab}_{\Gamma}\rightarrow G^{ab}_{\Gamma})\rightarrow 1
$$
of augmented pseudosimplicial groups yields, according to Theorem
\ref{t-10}(i), the exact sequence
\begin{align}\label{8}
0\rightarrow H^{\Gamma}_{2}(G)\rightarrow
\pi_{0}([P_{*},P_{*}]_{\Gamma})\overset{\tau'}{\rightarrow} G,
\end{align}
where the image of $\tau'$ is $[G,G]_{\Gamma}$. It is clear that
$P_{*}\rightarrow G$ is simplicially exact. For the calculation of
$\pi_{0}([P_{*},P_{*}]_{\Gamma})$ we will prove the equality
\begin{align}\label{9}
\Ker\partial^{1}_{0}\cap [P_{1},P_{1}]_{\Gamma}=
[P_{1},\Ker\partial^{1}_{0}]_{\Gamma}.
\end{align}
It suffices to show the inclusion $\Ker\partial^{1} _{0}\cap
[P_{1},P_{1}]_{\Gamma}\subset
[P_{1},\Ker\partial^{1}_{0}]_{\Gamma}$, since the converse
inclusion is obvious. First we will prove the equality
\begin{align}\label{10}
[P_{1},P_{1}]_{\Gamma}= [P_{1},\Ker\partial^{1}_{0}]_{\Gamma}\cdot
([P_{1},P_{1}]_{\Gamma}\cap s_{0}(P_{0})).
\end{align}
In effect, if $[x,y]_{\sigma}\in [P_{1},P_{1}]_{\Gamma}$, then for
$x= d\cdot s_{0}(c)$, $y= b\cdot s_{0}(a)$ with
$\partial^{1}_{0}(x)= c$, $\partial^{1}_{0}(y)= a$ and $b,d\in
\Ker\partial^{1}_{0}$, one has

\begin{align*}
&[x,y]_{\sigma}= [x,b\cdot s_{0}(a)]=
x\;^{\sigma}b\;^{\sigma}(s_{0}(a))x\;^{-1}\cdot s_{0}(a)^{-1}\cdot
b^{-1}= [x,b]_{\sigma}\cdot bx\;^{\sigma}(s_{0}(a))\cdot
x^{-1}\\
& \cdot s_{0}(a)^{-1}\cdot b^{-1}=[x,b]_{\sigma}\cdot
bds_{0}(c)\;^{\sigma}(s_{0}(a))s_{0}(c)^{-1}\cdot
d^{-1}s_{0}(a)^{-1}\cdot b^{-1}= [x,b]_{\sigma}\cdot
bd[s_{0}(c),s_{0}(a)]_{\sigma}\\
&\cdot s_{0}(a)d^{-1}s_{0}(a^{-1})\cdot b^{-1}=[x,b]_{\sigma}\cdot
z\cdot [s_{0}(c),s_{0}(a)]_{\sigma}\cdot bd\cdot
s_{0}(a)d^{-1}s_{0}(a)^{-1}b^{-1}= [x,b]_{\sigma}\cdot z\\ &\cdot
[s_{0}(c),s_{0}(a)]_{\sigma}\cdot b\cdot [d,s_{0}(a)]\cdot b^{-1}
\end{align*}
with $z\in [P_{1},\Ker\partial^{1}_{0}]_{\Gamma}$,
$b[d,s_{0}(a)]b^{-1}\in [P_{1},\Ker\partial^{1}_{0}]$ and
$[s_{0}(c),s_{0}(a)]_{\sigma}\in s_{0}(P_{0})$.

It follows that $[x,y]_{\Gamma}$ belongs to
$[P_{1},\Ker\partial^{1}_{0}]_{\Gamma}\cdot
([P_{1},P_{1}]_{\Gamma}\cap s_{0}(P_{0}))$ that proves the
required equality (\ref{10}).

Let $w\in [P_{1},P_{1}]_{\Gamma}\cap \Ker\partial^{1}_{0}$. Then
by (\ref{10}) $w\in [P_{1},\Ker\partial^{1}_{0}]\cdot
([P_{1},P_{1}]_{\Gamma}\cap s_{0}(P_{0}))$, that is $w= w'\cdot
x'$ with $w'\in [P_{1},\Ker\partial^{1}_{0}]_{\Gamma}$ and $x'\in
[P_{1},P_{1}]_{\Gamma}\cap s_{0}(P_{0})$. It follows that $x'=
w'^{-1}w$ belongs to $\Ker\partial^{1}_{0}\cap s_{0}(P_{0})=1$.
Whence $w= w'$ and the equality (\ref{9}) is proved.

Since $\partial^{1}_{1}(\Ker\partial^{1}_{0})= R$, using the
equality (\ref{9}) one gets
$$
\partial^{1}_{1}([P_{1},P_{1}]_{\Gamma}\cap
\Ker\partial^{1}_{0})=
\partial^{1}_{1}([P_{1},\Ker\partial^{1}_{0}]_{\Gamma})=
[P_{0},R]_{\Gamma}
$$
showing that $\pi_{0}([P_{*},P_{*}]_{\Gamma})= [P,P]_{\Gamma} /
[P,R]_{\Gamma}$. The required isomorphism of the theorem follows
now from the exact sequence (\ref{8}).
\end{proof}
\medskip

The result of Theorem 1.17(ii) will be called the equivariant Hopf
formula. In a forthcoming paper equivariant higher Hopf type
formulas will be established for the $\Gamma$-equivariant integral
homology $H_{*}^{\Gamma}(G)$.

\medskip

\begin{cor}
The group $R\cap [P,P]_{\Gamma} / [P,R]_{\Gamma}$ does not depend
on the projective presentation $P\rightarrow G$ of the
$\Gamma$-group $G$. A $\Gamma$-group has a universal central
$\Gamma$-equivariant extension if and only if it is
$\Gamma$-perfect.
\end{cor}

\medskip

\begin{thm}
Let $1\rightarrow N\rightarrow E\overset{\alpha}{\rightarrow}
G\rightarrow 1$ be a short exact sequence of $\Gamma$-groups such
that the $\Gamma$- homomorphism $\alpha$ has a $\Gamma$- excision
and $\tau:P\rightarrow E$ a projective presentation of the
$\Gamma$-group $E$.Then there is an exact sequence
$$
0\rightarrow U\rightarrow H^{\Gamma}_{2}(E)\rightarrow
H^{\Gamma}_{2}(G)\overset{\delta}{\rightarrow} N /
[E,N]_{\Gamma}\rightarrow H^{\Gamma}_{1}(E)\rightarrow
H^{\Gamma}_{1}(G)\rightarrow 0,
$$
where $U$ is the kernel of the $\Gamma$-homomorphism
$[P,S]_{\Gamma} / [P,R]_{\Gamma}\rightarrow [E,N]_{\Gamma}$, $R=
\Ker\tau$, $S= \Ker\alpha\tau$, induced by $\tau$.
\end{thm}
\begin{proof} Using the exact sequence (\ref{8}) the
$\Gamma$-homomorphism $\alpha$ induces the following commutative
diagram with exact rows and columns:
$$
\begin{matrix}
 & 0 & & 1 & & 1 & & 1 & & \\
 & \downarrow & & \downarrow & & \downarrow & & \downarrow & & \\
0\; \rightarrow &\Ker\gamma' & \rightarrow & \Ker\gamma &
\overset{\eta}{\rightarrow}& N & \rightarrow & N/\eta(\Ker\gamma)&
\rightarrow \; 0\\
 & \downarrow & & \downarrow & & \downarrow & & \downarrow & & \\
0\;  \rightarrow & H_2^{\Gamma}(E) & \rightarrow &
\pi_0([P_*,P_*]_{\Gamma}) & \rightarrow & E & \rightarrow &
H_1^{\Gamma}(E)&
\rightarrow \; 0\\

& \vcenter{\llap{$_{\gamma'}{}$}}\downarrow & &
\vcenter{\llap{$_{\gamma}{}$}}\downarrow & &
\downarrow\vcenter{\rlap{${}_{\alpha}$}} & &\downarrow && \\

0\; \rightarrow & H_2^{\Gamma}(G) & \rightarrow &
\pi_0([P^{G}_*,P^{G}_*]_{\Gamma}) & \rightarrow & G & \rightarrow
& H_1^{\Gamma}(G)&
\rightarrow \; 0\\
&  & & \downarrow & & \downarrow & & \downarrow & & \\
&  & & 1 & & 1 & & 0 & &
\end{matrix}
\;\;\;,
$$
where $P_{*}\rightarrow E$ and $P^{G}_{*}\rightarrow G$ are
$\mathcal{P}_{\mathcal{F}}$-projective resolutions of $E$ and $G$
induced by $\tau$ and $\alpha\tau$ respectively. Clearly
$\Ker\gamma= [P,S]_{\Gamma} / [P,R]_{\Gamma}$. It follows that the
image $(\eta(\Ker\gamma))$ of $\Ker\gamma$ is equal to
$[E,N]_{\Gamma}$ and $\Ker\eta$ is isomorphic to
$\Ker\gamma'$.Therefore this diagram yields the required exact
sequence, where the connecting homomorphism $\delta$ is defined in
a natural way.
\end{proof}

\medskip

Theorem 1.19 is a generalized equivariant version of the well
known Stallings-Stammbach exact sequence in integral homology
\cite{St}.

\medskip

\begin{thm}\label{t-20}
If $A$ is a $\Gamma$-equivariant $G$-module, then there is a
bijection
$$
E_{\Gamma}(G,A))\cong H^{2}_{\Gamma}(G,A).
$$
\end{thm}
\begin{proof} We will use the isomorphism
$H^{2}_{\Gamma}(G,A))\cong R^{1}_{\mathcal{F}}\Der_{\Gamma}(G,A))$
(see Theorem \ref{t-3}) and show the bijection
$E_{\Gamma}(G,A)\cong R^{1}_{\mathcal{F}}\Der(G,A)$. Take the free
cotriple resolution $F_{*}\rightarrow G$ of the $\Gamma$-group $G$
which is simplicially exact. Then $F_{1}\rightarrow F_{0}$ factors
through $F_{1}\rightarrow M\rightarrow F_{0}$, $l_{0}\tau_{1}=
\partial^{1}_{0}$, $l_{1}\tau_{1}=
\partial^{1}_{1}$, where $M$ is the simplicial kernel of $\tau$
and $\tau_{1}$ is surjective. If $f\in \Der_{\Gamma}(F_{1},A)$
such that $\sum f(-1)^{i}\partial^{2}_{i}= 0$, then there is
$f'\in \Der_{\Gamma}(M,A)$ such that $f'\tau_{1}= f$ and
$f'(\Delta)= 0$, $\Delta= \{(x,x), x\in F_{0}\}$. Denote by
$\widetilde{\Der}_{\Gamma}(M,A)$ the subgroup of
$\Der_{\Gamma}(M,A)$ consisting of $\Gamma$-derivations $f$ with
$f(\Delta)= 0$. Conversely, if $f'\in
\widetilde{\Der}_{\Gamma}(M,A)$, then $\underset{i}\sum
f'\tau_{1}(-1)^{i}\partial^{2}_{i}= 0$. It follows that it is
sufficient to establish a bijection with $\Coker\eta$, where
$\eta:\Der_{\Gamma}(F_{0},A)\rightarrow
\widetilde{\Der}_{\Gamma}(M,A)$, $\eta= \Der_{\Gamma}(l_{0},A)-
\Der_{\Gamma}(l_{1},A)$.

Define a map $\vartheta:E_{\Gamma}(G,A)\rightarrow
H^{2}_{\Gamma}(G,A)$ as follows. If $[E]\in E_{\Gamma}(G,A)$,
\linebreak$E: 0\rightarrow A\rightarrow B\rightarrow G\rightarrow
1 $, since $E$ is $\Gamma$-splitting, there is a commutative
diagram
$$
\begin{matrix}
M &\underset{l^{1}_{1}}{\overset{l^{1}_{0}}{\;\two\;}} & F_{0} &
\overset{\tau}{\rightarrow} & G \\
\vcenter{\llap{$_{f}{}$}}\downarrow & & \downarrow
\vcenter{\rlap{${}_{\varphi}$}} & & || \\
A &\overset{\alpha}{\rightarrow} & B &\overset{\beta}{\rightarrow}
& G  & ,
\end{matrix}
$$
where $F_{0}= F(G)$, $\varphi$  is a $\Gamma$-homomorphism and
$f(x)= \varphi l_{0}(x)\cdot \varphi l_{1}^{-1}(x)$, $x\in M$. It
is easily seen that $f$ is a $\Gamma$-derivation such that
$f(\Delta)= 0$ and define $\vartheta([E])= [f]$. Conversely, if
$[f]\in H^{2}_{\Gamma}(G,A)$, then take the semidirect product
$A\rtimes G$ and introduce a relation
$$
(a,x)\sim (a',x')\Leftrightarrow \tau(x)=\tau(x')
$$
and $a\cdot f(x,x')= a'$. It is easy to check that this relation
is a congruence and let $C$ be the quotient $A\rtimes G / \sim$
which is a $\Gamma$-group. One gets a commutative diagram
$$
\begin{matrix}
M &\underset{l^{1}_{1}}{\overset{l^{1}_{0}}{\;\two}} & F(G) &
\overset{\tau}{\rightarrow} & G \\
\vcenter{\llap{$_{f}{}$}}\downarrow & & \downarrow
\vcenter{\rlap{${}_{\psi}$}} & || \\
 A &\overset{\sigma}{\rightarrow} & C &\overset{\mu}{\rightarrow}
 & G & ,
\end{matrix}
$$
where  $\sigma (a)= [(a,x)]$, $\mu([(a,x)])= \tau(x)$, $\psi(x)=
[(o,x)]$. The extension $E : 0\rightarrow
A\overset{\sigma}{\rightarrow} C\overset{\mu}{\rightarrow}
G\rightarrow 1$ is a $\Gamma$-equivariant extension of $G$ by $A$,
the splitting $\Gamma$-map is given by $\gamma(g)= \psi(0,|g|)$,
$g\in G$. Define $\vartheta': H^{2}_{\Gamma}(G,A)\rightarrow
E_{\Gamma}(G,A)$ by $\vartheta'([f])= [E]$. It is standard to show
that $\vartheta$ and $\vartheta'$ are well defined and inverse to
each other.
\end{proof}

\medskip

Note that Theorem \ref{t-20} could also be proved using the
corresponding factor set theory for $\Gamma$-groups and the
bijection $\vartheta$ is in fact an isomorphism with respect to
the "Baer sum" which could be introduced on $E_{\Gamma}(G,A)$. The
description of higher $\Gamma$-equivariant group cohomology
$H^{n+1}_{\Gamma}(G,A)$, $n\geq 2$, by extensions is also
realizable using $n$-fold $\Gamma$-equivariant extensions of $G$
by $A$, that means extensions of the form
$$
E: 0\rightarrow A\rightarrow X_{1}\rightarrow X_{2}\rightarrow
\cdots \rightarrow X_{n}\rightarrow G\rightarrow 1,
$$
where $0\rightarrow A\rightarrow X_{1}\rightarrow
Im\alpha_{1}\rightarrow 0$, $0\rightarrow
Im\alpha_{i-1}\rightarrow X_{i}\rightarrow Im\alpha_{i}\rightarrow
0$, for \linebreak$1\leq i \leq n-1$, are proper short exact
sequences of $G\rtimes \Gamma$-modules and \linebreak$0\rightarrow
Im\alpha_{n-1}\rightarrow X_{n}\rightarrow G\rightarrow 1$ is a
$\Gamma$-equivariant extension of $G$ by $Im\alpha_{n-1}$, and by
introducing the $\Gamma$-equivariant characteristic class
$\chi(E)$ of a $\Gamma$-equivariant extension $E$ of $G$ by $A$,
which means by constructing for any $\Gamma$-equivariant extension
$E$ of $G$ by $A$:
$$
E: 0\rightarrow A\overset{\alpha}{\rightarrow}
B\overset{\beta}{\rightarrow} G\rightarrow 1
$$
the exact sequence of $G\rtimes \Gamma$-modules
$$
\chi(E)= 0\rightarrow A\overset{\alpha'}{\rightarrow}
F(B')/L\overset{\beta'}{\rightarrow}
\mathbb{Z}(G)\overset{\epsilon}{\rightarrow} \mathbb{Z}\rightarrow
0,
$$
where $F(B')$ is the relatively free $G\rtimes \Gamma$-module
generated by the $\Gamma$-set $B'= \{[b], b\in B, b\neq 0\}$ with
$[0]= 0$, $^{\sigma}(x[b])=\; ^{\sigma}x[^{\sigma}b]$ for $x\in
G$, $\sigma\in \Gamma$, $b\in B'$; $L$ is its $G\rtimes
\Gamma$-submodule generated by the elements
$[b_{1}+b_{2}]-\beta(b_{1})[b_{2}]- [b_{1}]$, where $b_{1},
b_{2}\in B$, and the $G\rtimes \Gamma$-homomorphisms $\alpha',
\beta'$ are induced in a natural way by $\alpha$ and $\beta$
respectively.

 Using the cochain description (see Section 1)
of the $\Gamma$-equivariant cohomology of groups
$H^{*}_{\Gamma}(G,A)$ the cup product can be defined, since the
tensor product \cite{AW} of $\Gamma$-cochains is again a
$\Gamma$-cochain. Therefore there is a cup product
\begin{align*}
H^{p}_{\Gamma}(G,A)\otimes H^{q}_{\Gamma}(G,B)\rightarrow
H^{p+q}_{\Gamma}(G,A\otimes B)
\end{align*}
for $p,q\geq 1$, endowing on $H^{*}_{\Gamma}(G,A)$ a structure of
$H^{*}_{\Gamma}(G)$-module.

For any short exact sequence of $G\rtimes \Gamma$-modules

$$
 0\rightarrow A'\overset{\alpha}{\rightarrow}
A\overset{\beta}{\rightarrow} A''\rightarrow 0
$$
such that $\beta$ is $\Gamma$-splitting the sequences
$$
0\rightarrow P\otimes_{G\rtimes \Gamma}A'\rightarrow
P\otimes_{G\rtimes \Gamma}A\rightarrow P\otimes_{G\rtimes
\Gamma}A''\rightarrow 0
$$
and
$$
0\rightarrow Hom_{G\rtimes \Gamma}(P,A')\rightarrow Hom_{G\rtimes
\Gamma}(P,A)\rightarrow Hom_{G\rtimes \Gamma}(P,A'')\rightarrow 0
$$
are exact for any relatively projective $G\rtimes \Gamma$-module
$P$ implying exact $\Gamma$-equivariant homology and cohomology
sequences respectively.

Let $G$ be finite and let us consider the homomorphism $N_{G}:
A\rightarrow A$, $N_{G}(a)=\underset{s\in G}{\Sigma}^{s}a$, where
$A$ is a $G\rtimes \Gamma$-module. Assume that $\Gamma$ acts
trivially on $N_{G}(A)$. Therefore one has the inclusion
$N_{G}(A)\subset A^{G\rtimes \Gamma}$ and $N_{G}$ induces a
homomorphism $N_{G}^{*}: H_{0}^{\Gamma}(G,A)\rightarrow
H_{\Gamma}^{0}(G,A)$. For this we have only to show that
$N_{G}(^{\sigma}a-a)= 0$. In effect $N_{G}(a)= ^{\sigma}N_{G}(a)=
^{\sigma}(\underset{s\in G}{\Sigma}^{s}a)= \underset{s\in
G}{\Sigma}^{\sigma}(^{s}a)= \underset{s\in
G}{\Sigma}^{^{\sigma}s}(^{\sigma}a)= N_{G}(^{\sigma}a)$. Under the
afore given assumption we can define $\Gamma$-equivariant Tate
cohomology groups $\check{H}_{\Gamma}^{n}(G,A)$, $n\in
\mathbb{Z}$, by setting
$$
\check{H}_{\Gamma}^{n}(G,A)= H_{\Gamma}^{n}(G,A), n\geq 1,
$$
$$
\check{H}_{\Gamma}^{0}(G,A)= \Ker N_{G}^{*},\;
\check{H}_{\Gamma}^{-1}(G,A)= \Coker N_{G}^{*},
$$
$$
\check{H}_{\Gamma}^{-n}(G,A)= H_{n-1}^{\Gamma}(G,A), n\geq 2.
$$

\begin{prop}
 For any short exact sequence of $G\rtimes \Gamma$-modules
$$
 E: 0\rightarrow A'\overset{\alpha}{\rightarrow}
A\overset{\beta}{\rightarrow} A''\rightarrow 0
$$
such that $\beta$ is $\Gamma$-splitting there is a long exact
sequence of $\Gamma$-equivariant Tate cohomology groups
$$
\cdots \rightarrow \hat{H}^{n-1}_{\Gamma}(G,A'')\rightarrow
\hat{H}^{n}_{\Gamma}(G,A')\rightarrow
\hat{H}^{n}_{\Gamma}(G,A)\rightarrow
\hat{H}^{n}_{\Gamma}(G,A'')\rightarrow
\hat{H}^{n+1}_{\Gamma}(G,A')\rightarrow\cdots
$$
\end{prop}
\begin{proof}
By using for $E$ the exact $\Gamma$-equivariant homology and
cohomology sequences the proof is similar to the classical case
(see [1], Chapter IV, Theorem 6.1).
\end{proof}

 Note that it would be interesting to construct the
$\Gamma$-equivariant versions of Farell cohomology theory of
groups \cite{BrK} and Vogel (co)homology theory of groups
\cite{Vgl} generalizing the above defined $\Gamma$-equivariant
Tate cohomology theory of groups.

\

\

\section{Relationship with equivariant cohomology of topological
spaces}

\

Let $X$ be a topological space. If a group $G$ acts on $X$, then
this action induces an action of $G$ on the singular complex
$S(X)$ of $X$ given by $gf$, $f:\Delta_{n}\rightarrow X$, $f\in
S(X)$, making $S(X)$ a chain complex of $G$-modules, where
$g:X\rightarrow X$ is the homeomorphism of $X$ induced by the
action of $g\in G$.

Throughout out this section $X$ is a $G$-space the group $G$
acting on $X$ properly. That means each point $x$ of $X$ belongs
to some proper open subset of $X$. Recall that an open subset $U$
of $X$ is called proper with respect to the action of $G$, if
$^{g}U\cap U$ is empty for all elements $g\neq 1$ of $G$
\cite{Ml}.

We will assume also that a separate group $\Gamma$ acts on $G$ and
$X$ such that the following condition holds:
\begin{align}\label{13}
^{\sigma}(^{g}x)= ^{^{\sigma}g}(^{\sigma}x),
\end{align}
for $x\in X$, $g\in G$, $\sigma\in \Gamma$.

For example if $X$ is a $G$-space, take $\Gamma= G$ with the
actions of $\Gamma$ on $G$ by conjugation and on $X$ as $G$ is
acting.

Then the augmented singular complex $S(X)\rightarrow Z$:
\begin{align}\label{14}
\cdots \rightarrow S_{n}(X)\rightarrow S_{n-1}(X)\rightarrow
\cdots \rightarrow S_{1}(X)\rightarrow S_{0}(X)\rightarrow
\mathbb{Z}\rightarrow 0
\end{align}
is a chain complex of $G\rtimes \Gamma$-modules, where the groups
$G$ and $\Gamma$ act trivially on $\mathbb{Z}$.

It will be said that the topological space $X$ has the property
$(c)$, if the singular complex (\ref{14}) is exact and any induced
short exact sequence
$$
0\rightarrow \Ker\partial_{k}\rightarrow S_{k}(X)\rightarrow
Im\partial_{k}\rightarrow 0
$$
is $\Gamma$-splitting for $k\geq 0$.

For instance $X$ satisfies the condition $(c)$ if either $X$ is
acyclic and $\Gamma$ acts trivially on $X$ or $X$ is
$\Gamma$-contractible, that is the identity map
$1_{X}:X\rightarrow X$ is $\Gamma$-homotopic to a constant map
$f_{0}:X\rightarrow x_{0}\in X$.

\medskip

\begin{thm}\label{t-23} If a topological space $X$ satisfies the condition
$(c)$, then there is an isomorphism
$$
H^{n}_{\Gamma}(G,A)\cong H^{n}_{\Gamma}(X / G ,A)
$$
for $n\geq 0$, where $A$ is an abelian group on which $G$ and
$\Gamma$ act trivially and $H^{*}_{\Gamma}(X / G,A)$ is the
equivariant cohomology of topological spaces \cite{BrG}.
\end{thm}
\begin{proof} Since $G$ acts properly on $X$, by Lemma
11.2 (\cite{Ml},Chapter IV) the sequence (\ref{14}) is a chain
complex of free $G$-modules. Therefore each $S_{n}(X)$, $n\geq 0$,
is a $G\rtimes \Gamma$-module which is free as $G$-module and its
basis consisting of singular $n$-th simplexes is a
$\Gamma$-subset. It follows that (\ref{14}) is a relatively free
$G\rtimes \Gamma$-resolution of $\mathbb{Z}$.

By Proposition 11.4 (\cite{Ml}, Chapter IV) the canonical map
$p:X\rightarrow X / G$ induces an isomorphism
$$
p^{*}: \Hom_{Z}(S(X / G),A)\cong \Hom_{G}(S(X),A)
$$
of chain complexes. Notice that the group $\Gamma$ acts naturally
on $X / G$ and the map $p$ is a $\Gamma$-map. Indeed, the action
of $\Gamma$ given by $^{\sigma}([x])= [^{\sigma}x]$, $x\in X$,
$\sigma\in \Gamma$, is well defined, thanks to the equality
(\ref{13}); if $^{g}x= y$ for some $g\in G$, then $^{\sigma}y=\;
^{\sigma}(^{g}x)= \;^{^{\sigma}g}(^{\sigma}x)$ for any $\sigma\in
\Gamma$. It is obvious that under so defined action of $\Gamma$
the map $p$ is a $\Gamma$-map. This implies that the isomorphism
$p^{*}$ induces an isomorphism
$$
\Hom_{\Gamma}(S(X / G),A)\cong \Hom_{G\rtimes \Gamma}(S(X),A)
$$
of cochain complexes giving the required isomorphism of the
equivariant group cohomology of the space $X / G$ with the
equivariant group cohomology of the group $G$ with coefficients in
$A$ with trivial actions of $G$ and $\Gamma$.
\end{proof}

\medskip

Property $(c)$ holds whenever $\Gamma$ acts trivially on the group
$G$ and on the acyclic space $X$. Therefore Theorem \ref{t-23} is
an equivariant version of the classical Theorem 11.5
(\cite{Ml},Chapter IV).

\

\

\section{Applications to algebraic K-theory}

\

Let $A$ be a unital ring and $I$ its ideal. The group $E_{n}(A,I)$
is the normal subgroup of the group $E_{n}(A)$ of elementary
$n$-matrices generated by $I$-elementary $n$-matrices
$\epsilon_{n}$ of the form $\epsilon_{n}= I_{n}+ ae_{ij}$ with
$a\in I$ and $i\neq j$ (see \cite{Ba}). The group $E(A,I)$ is
defined as $\underset{\underset{n}{\rightarrow}}{\lim}
E_{n}(A,I)$. It is known \cite{Ba} that
$$
E_{n}(A,I)= [E_{n}(A), E_{n}(A,I)]
$$
for $n\geq 3$.

It follows that the group $E_{n}(A,I)$, $n\geq 3$, which is not
perfect in general, is a $E_{n}(A)$-perfect group, the group
$E_{n}(A)$ acting on $E_{n}(A,I)$ by conjugation. Clearly the same
is true for the relative elementary group $E(A,I)$, that is
$E(A,I)$ is a $E(A)$-perfect group.

Now let $A$ be a ring not necessarily with identity. Denote by
$A^{+}$ the unital ring given by $A^{+}= \{(a,n), a\in A, n\in
\mathbb{Z}\}$ with usual sum and product
$$
(a,n)\cdot (a',n')= (aa'+ na'+ n'a, nn').
$$
One has a short split exact sequence of rings
\begin{align}\label{15}
0\rightarrow A\overset{\sigma}{\rightarrow}
A^{+}\overset{\tau}{\rightarrow} \mathbb{Z}\rightarrow 0,
\end{align}
where $\sigma(a)= (a,0)$, $\tau(a,n)= n$ and with splitting map
$\gamma:Z\rightarrow A^{+}$, $\gamma(n)=(0,n)$.

By definition $E(A)= \Ker E(\tau)$, $St(A)= \Ker St(\tau)$ and
$K_{2}(A)= \Ker K_{2}(\tau)$. Clearly $E(A)= E(A^{+},A)$. Whence
we have the following short exact sequence of short exact
sequences induced by (\ref{15}):
\begin{align*}
0\rightarrow (0\rightarrow K_{2}(A)\rightarrow
St(A)\overset{\beta}{\rightarrow} E(A)\rightarrow 1) \rightarrow
(0\rightarrow K_{2}(A^{+})\rightarrow
St(A^{+})\overset{\beta^{+}}{\rightarrow} E(A^{+})\rightarrow 1)\\
\rightarrow (0\rightarrow K_{2}(\mathbb{Z})\rightarrow
St(\mathbb{Z})\rightarrow E(\mathbb{Z})\rightarrow 1) \rightarrow
0.
\end{align*}

One gets an action of $St(\mathbb{Z})$ on $St(A)$ by conjugation
using the splitting map $St(\gamma)$. The group $E(A)$ and the
general linear group $GL(A)$ also become $St(\mathbb{Z})$-groups
by conjugation via the map $\beta^{+}\cdot St(\gamma)$. Clearly
$\beta$ is a $St(\mathbb{Z})$-homomorphism and $St(\mathbb{Z})$
acts trivially on $K_{2}(A)$. Therefore the central extension of
the group $E(A)$
\begin{align}\label{16}
0\rightarrow K_{2}(A)\rightarrow St(A)
\overset{\beta}{\rightarrow} E(A)\rightarrow 1
\end{align}
is a sequence of $St(\mathbb{Z})$-groups.

There is a presentation of the Steinberg group $St(A)$ as a
$St(\mathbb{Z})$-group as follows \cite{Sw1}. The generators
$x^{a}_{ij}$ for $i,j\geq 1$, $i\neq j$, $a\in A$, satisfy the
relations
\begin{enumerate}
\item[(1)] $x^{a}_{ij}x^{b}_{ij}= x^{a+b}_{ij}$.
\item[(2)] $[x^{a}_{ij},x^{b}_{km}] = 1, j\neq k, i\neq m$.
\item[(3)] $[x^{a}_{ij},x^{b}_{jk}] = x^{ab}_{ik}$.
\item[(4)] $x^{z}_{ij}x^{a}_{ij}(x^{z}_{ij})^{-1} = x^{a}_{ij}, z\in
Z$.
\item[(5)] $x^{z}_{ij}x^{b}_{km}(x^{z}_{ij})^{-1} = x^{b}_{km}, z\in Z, i\neq m, j\neq
k$.
\item[(6)] $x^{z}_{ij}x^{b}_{jk}(x^{z}_{ij})^{-1} = x^{zb}_{ik}x^{b}_{jk}, z\in Z, i\neq
k$.
\item[(7)] $x^{z}_{ij}x^{b}_{ki}(x^{z}_{ij})^{-1} = x^{-zb}_{kj}x^{b}_{ki}, z\in Z, j\neq k$.
\end{enumerate}

\medskip

\begin{thm}\label{t-25}
The sequence (\ref{16}) is a universal central
$St(\mathbb{Z})$-equivariant extension of $E(A)$.
\end{thm}
\begin{proof} The map $\kappa^{+}:
\{e^{a^{+}}_{ij}\}\rightarrow St(A^{+})$ of elementary matrices
given by $\kappa^{+}{(e^{a^{+}}}_{ij})= x^{a^{+}}_{ij}$, $a^{+}\in
A^{+}$, $i\neq j$, induces a $St(\mathbb{Z})$- map $\kappa:
E(A)\rightarrow St(A)$ such that $\beta \kappa= 1_{E(A)}$ by using
the proof of a similar fact in the case of Corollary 15 and taking
into account the action of $St(\mathbb{Z})$ on the generators
$e_{ij}^{a}$. We conclude that (\ref{16}) is a central
$St(\mathbb{Z})$-equivariant extension of $E(A)$.

From the above condition (6) on generators one has
$x^{1}_{ij}x^{b}_{jk}(x^{1}_{ij})^{-1} = x^{b}_{ik}x^{b}_{jk}$,
$i\neq k$, $b\in A$. It follows that $[St(A),St(\mathbb{Z})]=
St(A)$ showing the groups $St(A)$ and $E(A)$ are
$St(\mathbb{Z})$-perfect groups. Following Theorem \ref{t-16}, it
remains to show that any central $St(\mathbb{Z})$-equivariant
extension of $St(A)$ splits.

Let
\begin{align}\label{17}
o\rightarrow C\rightarrow Y\overset{\theta}{\rightarrow}
St(A)\rightarrow 1
\end{align}
be a central $St(\mathbb{Z})$-equivariant extension of $St(A)$
with splitting $St(\mathbb{Z})$-map $\bar{\gamma}:
St(A)\rightarrow Y$. Consider the following exact sequence
\begin{align}\label{18}
0\rightarrow C\rightarrow Y\rtimes St(\mathbb{Z})\rightarrow
St(A)\rtimes St(\mathbb{Z})\rightarrow 1
\end{align}
induced by (\ref{17}) and by the given action of $St(\mathbb{Z})$
on this sequence. Clearly $St(A^{+}) \cong St(A)\rtimes
St(\mathbb{Z})$ and (\ref{18}) is a central extension of
$St(A^{+})$. Therefore (\ref{18}) splits, since $\beta^{+}:
St(A^{+})\rightarrow E(A^{+})$ is a universal central extension of
$E(A^{+})$ . Obviously $St(A^{+})$ is generated by $x^{a}_{ij}$
and $x^{z}_{kl}$, $a\in A$, $z\in \mathbb{Z}$.

Following the proof of Theorem 5.10 \cite{Mi}, the section for
(\ref{18}) can be constructed as follows. For $i\neq j$ chose an
index $h$ distinct from $i$ and $j$. Take the elements $y=
\bar{\gamma}(x^{1}_{ih})$ and $y'= \bar{\gamma}(x^{a}_{hj})$. Then
the needed section is defined by $x^{a}_{ij}\mapsto [y,y'] =
s^{a}_{ij}$, $x^{z}_{jk}\mapsto x^{z}_{jk}$. In \cite{Mi} it is
shown that this section does not depend on $h$ and the elements
$s^{a}_{ij}, x^{z}_{jk}$ satisfy all the Steinberg relations. This
implies that the elements $s^{a}_{ij}, x^{z}_{jk}$ satisfy the
relations (1)-(7) of the $St(\mathbb{Z})$-presentation of $St(A)$.
It follows that by sending $x^{a}_{ij}$ to $s^{a}_{ij}$, $a\in A$,
$i\neq j$, this map gives rise to the required  splitting
$St(\mathbb{Z})$-homomorphism $s:St(A)\rightarrow Y$.
\end{proof}

\medskip

\begin{cor}
There is an isomorphism $$ K_{2}(A)\cong
H^{St(\mathbb{Z})}_{2}(E(A))\quad \text{and}\quad
H^{St(\mathbb{Z})}_{2}(St(A))= 0$$ for any ring $A$.
\end{cor}
\begin{proof} The proof follows from Theorems \ref{t-17}, \ref{t-25} and \ref{t-20}.
\end{proof}

Finally we provide the construction of an alternative equivariant
algebraic K-theory $K^{\Gamma}_{*}$ by using $\Gamma$-equivariant
commutators.

Let $\Gamma_{i}(G)$ , $i\geq 0$, be the lower $\Gamma$-equivariant
central series of a $\Gamma$-group $G$ \cite{CeIn}, where
$\Gamma_{0}(G)= G$, $\Gamma_{1}(G)= [G,G]_{\Gamma}$ and
$\Gamma_{i+1}(G)= [G,\Gamma_{i}(G)]_{\Gamma}$, $i\geq 0$. First we
give the equivariant version of the $Z$-completion functor
$Z_{\infty}:\mathcal{G}\rightarrow \mathcal{G}$ defined on the
category of groups \cite{BK}, by setting
$$
Z^{\Gamma}_{\infty}(G)=
\underset{\underset{i}{\leftarrow}}{\lim}\;G / \Gamma_{i}(G),
$$
where $\{\Gamma_{i}(G)\}$ is the $\Gamma$-equivariant lower
central series of the $\Gamma$-group $G$. We obtain a covariant
functor $Z^{\Gamma}_{\infty}:\mathcal{G}_{\Gamma}\rightarrow
\mathcal{G}_{\Gamma}$.

Let $A$ be a ring and $\Gamma$ a group acting on the general
linear group $GL(A)$. Define the $\Gamma$-equivariant algebraic
K-functors by
\begin{align}\label{19}
K^{\Gamma}_{n}(A)=
L^{{\mathcal{P}}_{\mathcal{F}}}_{n-1}Z^{\Gamma}_{\infty}(GL(A)),
n\geq 1,
\end{align}
where $\mathcal{P}_{\mathcal{F}}$ is the projective class induced
by the free cotriple $\mathcal{F}$ in the category
$\mathcal{G}_{\Gamma}$ of $\Gamma$-groups. This definition could
actually be considered as an equivariant version of Quillen's
algebraic K-theory, since in the case of the trivial action of
$\Gamma$ on $GL(A)$ it is proved that the left derived functors of
the functor $Z\infty$ with respect to the projective class induced
by the free cotriple in the category of groups are isomorphic to
Quillen's $K$-groups up to dimension shift \cite{Ke, InH2}. It
would be interesting to establish the relationship of the afore
defined equivariant algebraic $K$-theory (\ref{19}) with the
equivariant algebraic {K}-theory given in \cite{FHM}.

\

\

\begin{center}
\bf Acknowledgments
\end{center}

I would like to thank the referee for careful examination of the
first draft of the manuscript, for useful comments and
suggestions. The author was partially supported by INTAS grant No
566, FNRS grant No 7GEPJ06551301 and  NATO linkage grant PST.CLG.
979167

\

\


\begin{thebibliography}{[*****]}

\bibitem{AW} M.F.Atiyah and C.T.C.Wall, {\em Cohomology of groups},
 Algebraic Number Theory, J.W.S.Cassels and A.Frohlich editors,
London and New York, Academic Press, 1967.

\bibitem{Ba} H.Bass, {\em Algebraic K-theory}, W.A.Benjamin,
Inc., New York, 1968.

\bibitem{BK} A.K.Bousfield and D.M.Kan,{\em Homotopy with respect to a
ring}, Algebraic Topology, Proc.Symp.in Pure Mathematics {\bf 22}
(1971).

\bibitem{BrG} G.Bredon, {\em Equivariant cohomology theories},
Lecture Notes in Math., Springer-Verlag, {\bf 34} (1967).

\bibitem{BrK} K.S.Brown, {\em Cohomology of groups}, Graduated Texts in Math.,
Springer-Verlag, 1982.

\bibitem{Ca} G.Carlsson, {\em Equivariant Stable Homotopy Theory and Related
Areas}, Proc. Workshop at Stanford University 2000, Carlsson
editor, Homology, Homotopy and Applications {\bf 3}(2), (2001).

\bibitem{CGO1} A.M.Cegarra, J.M.Garcia-Calcines and
J.A.Ortega, {\em On graded categorical groups and equivariant
group extensions}, Canadian J. of Math., {\bf 54} (5) (2002),
970-997.

\bibitem{CGO2} A.M.Cegarra, J.M.Garcia-Calcines and
J.A.Ortega, {\em Cohomology of groups with operators}, Homology,
Homotopy and Applications {\bf 4} (1) (2002), 1-23.

\bibitem{CeGa1} A.M.Cegarra and A.R.Garzon,
{\em Equivariant group cohomology and Brauer group}, Bull. Belg.
Math. Soc. Simon Stevin {\bf 10} (3) (2003), 451-459.

\bibitem{CeGa2} A.M.Cegarra and A.R.Garzon,
 {\em Some algebraic applications of graded categorical groups
 theory}, Theory Appl. Categ. {\bf 11} (10) (2003), 215-251.

\bibitem{CeIn} A.M.Cegarra and H.Inassaridze, {\em Homology of groups with
operators}, Intern. Math. J., {\bf 5}, (1) (2004), 29-48.

\bibitem{El} S.Eilenberg, {\em Foundations of relative homological
algebra}, Memoirs Amer. Math. Soc. {\bf 55} (1965), 1-39.

\bibitem{FHM} Z.Fiedorowicz, H.Hauschild and J.P.May,
{\em Equivariant algebraic K-theory}, Lecture Notes in Math.,
Springer-Verlag {\bf 967} (1982), 23-80.

\bibitem{Fr1} A.Frabetti, {\em Leibniz homology of dialgebras of
matrices}, JPAA {\bf 129} (1998), 123-141.

\bibitem{Fr2} A.Frabetti, {\em Dialgebra (co)homology with coefficients},
1999, preprint.

\bibitem{FrhWa} A.Frohlich and C.T.C.Wall, {\em Graded monoidal
categories}, {\bf 28} (1974), 229-285.

\bibitem{Ge} S.M.Gersten, {\em $K_{3}$ of a ring is $H_{3}$ of the Steinberg
Group}, Proc. Amer. Math. Soc., {\bf 37} (2) (1973), 366-368 .

\bibitem{Go} F.Goichot, {\em Homologie de Tate-Vogel
equivariante}, JPAA {\bf 82} (1992), 39-64.

\bibitem{InH1} H.Inassaridze, {\em Homotopy of pseudosimplicial groups,
non-abelian derived functors and algebraic K-theory} Math. USSR
Sbornik {\bf 98}, No3 (1975), 339-362 (in Russian).

\bibitem{InH2} H.Inassaridze, {\em Non-Abelian Homological Algebra and its
Applications}, Kluwer Academic Publishers, 1997.

\bibitem{InN} N.Inassaridze, {\em Non-abelian tensor products and non-abelian
homology of groups}, JPAA {\bf 112} (1996), 191-205.

\bibitem{InHInN} H.Inassaridze and N.Inassaridze, {\em Non-abelian homology of
groups}, K-Theory, {\bf 18} (1999), 1-17.

\bibitem{Ke} F.Keune, {\em Derived functors and
algebraic K-theory}, Algebraic K-theory I, Battelle Institute
Conference 1972, Lecture notes in Math., Springer-Verlag, {\bf
341} (1973), 166-176.

\bibitem{Ku} A.Kuku, {\em Equivariant K-theory and the cohomology of profinite
groups}, Algebraic K-theory, Number Theory, Geometry and Analysis,
Lecture Notes in Math., Springer-Verlag, {\bf 1046} (1984),
235-244.

\bibitem{Lo1} J.-L.Loday, {\em Dialgebras}, Prepublication 1999/14,
Institut de Recherche Mathematique, Universite Louis Pasteur,
Strasbourg, 1999.

\bibitem{Lo2} J.-L.Loday, {\em Cohomologie et groupe de Steinberg
relatifs}, J. Algebra, {\bf 54} (1978), 178-202.

\bibitem{Ml} S.MacLane, {\em Homology}, Springer-Verlag, 1963.

\bibitem{Mi} J.Milnor, {\em Introduction to algebraic
K-theory}, Princeton University Press, 1971.

\bibitem{Ph} N.C.Philips, {\em Equivariant K-theory and freeness
of group actions on C*-algebras}, Lecture Notes in Math.,
Springer-Verlag, {\bf 1274}, 1987.

\bibitem{Qu} D.Quillen, {\em Spectral sequences of a double semi-simplicial
group}, Topology {\bf 5} (1966),155-157.

\bibitem{St} U.Stammbach, {\em Homology in group theory},
Lecture Notes in Math., Springer-Verlag, {\bf 359}, 1973.

\bibitem{Sw1} R.Swan, {\em Nonabelian homological algebra and
K-theory}, Proc. Symposia in Pure Mathematics {\bf 17} (1970),
88-123.

\bibitem{Sw2} R.Swan, {\em Some relations between higher
K-functors}, J. Algebra {\bf 21}, No1 (1972), 113-136.

\bibitem{TiVo} M.Tierney and W.Vogel, {\em Simplicial resolutions and derived
functors}, Math. Zeit. {\bf 111}(1969), 1-14.

\bibitem{Vgl} P.Vogel, {\em Talks given at Paris and Rennes
1983-84}.

\bibitem{Wh} J.H.C.Whitehead, {\em On group extensions with
operators}, Quart. J. Math. Oxford {\bf (2)}, 1(1950), 219-228.

\end{thebibliography}
\end{document}